\documentclass[journal,web]{ieeecolor}
\usepackage{IEEETAC}
\usepackage[compress]{cite}
\usepackage{amsmath,amssymb,amsfonts}
\usepackage{algorithmic}
\usepackage{graphicx}
\usepackage{textcomp}
\usepackage{comment} 
\usepackage{multirow}
\usepackage{tabularx}
\usepackage{float} 
\usepackage{setspace}

\DeclareMathOperator*{\argmax}{arg\,max}

\def\BibTeX{{\rm B\kern-.05em{\sc i\kern-.025em b}\kern-.08em
    T\kern-.1667em\lower.7ex\hbox{E}\kern-.125emX}}
\begin{document}

\begin{titlepage}
\noindent IEEE Copyright Notice \\

\noindent © 2020 IEEE.  Personal use of this material is permitted.  Permission from IEEE must be obtained for all other uses, in any current or future media, including reprinting/republishing this material for advertising or promotional purposes, creating new collective works, for resale or redistribution to servers or lists, or reuse of any copyrighted component of this work in other works.\\

\noindent Citation information: DOI 10.1109/TAC.2020.3027804, IEEE Transactions on Automatic Control.
\end{titlepage}

\title{Fast Algorithm for Fuel-Optimal Impulsive Control of Linear Systems with Time-Varying Cost}
\author{Adam W. Koenig and Simone D'Amico
\thanks{This work was supported by a NASA Office of the Chief Technologists Space Technology Research Fellowship (NSTRF), NASA Grant \#NNX15AP70H and the U.S. Air Force Research Laboratory's Control, Navigation, and Guidance for Autonomous Spacecraft (CoNGAS) contract FA9453-16-C-0029.}
\thanks{Adam W. Koenig (e-mail: awkoenig@stanford.edu) and Simone D'Amico (e-mail: damicos@stanford.edu) are with the Department of Aeronautics and Astronautics, Stanford University, Stanford, CA 94305 USA.}}

\maketitle


\begin{abstract}
This paper presents a new fast and robust algorithm that provides fuel-optimal impulsive control input sequences that drive a linear time-variant system to a desired state at a specified time.
This algorithm is applicable to a broad class of problems where the cost is expressed as a time-varying norm-like function of the control input, enabling inclusion of complex operational constraints in the control planning problem.
First, it is shown that the reachable sets for this problem have identical properties to those in prior works using constant cost functions, enabling use of existing algorithms in conjunction with newly derived contact and support functions.
By reformulating the optimal control problem as a semi-infinite convex program, it is also demonstrated that the time-invariant component of the commonly studied primer vector is an outward normal vector to the reachable set at the target state.
Using this formulation, a fast and robust algorithm that provides globally optimal impulsive control input sequences is proposed.
The algorithm iteratively refines estimates of an outward normal vector to the reachable set at the target state and a minimal set of control input times until the optimality criteria are satisfied to within a user-specified tolerance.
Next, optimal control inputs are computed by solving a quadratic program.
The algorithm is validated through simulations of challenging example problems based on the recently proposed Miniaturized Distributed Occulter/Telescope small satellite mission, which demonstrate that the proposed algorithm converges several times faster than comparable algorithms in literature.
\end{abstract}



\section{Introduction}\label{sec:intro}

The fuel-optimal impulsive control problem for linear time-variant dynamical systems with fixed end times and states has received a great deal of attention in literature.
In this paper, the term ``fuel-optimal" means that the cost metric is expressed as the integral of a norm-like function of the control input and has no dependence on the state.
This is consistent with resource costs in a number of practical applications such as a propellant expenditure for a spacecraft thruster executing small maneuvers or energy expenditure in an electric motor in a drone.
The additional feature of impulsive control is that the magnitude of the control input is not constrained.
However, a wide range of continuous control problems can be approximated as impulsive provided that the durations of the time intervals over which control input is applied are small.
Because of these properties, similar problems have been studied in a wide range of fields including engineering \cite{Gaias2015diffdrag1,Chernick2018}, epidemiology \cite{Verriest2005}, and finance \cite{Bouchard2011}.
Indeed, this class of problem has been studied for over fifty years in the context of spacecraft rendezvous and formation-flying \cite{Prussing1969}.
The space community's interest in these problems is motivated by the fact that spacecraft propellant is limited and cannot be replenished after launch.
It follows that improving the efficiency of maneuver planning algorithms can significantly extend mission lifetimes.
Additionally, the dynamics of the space environment are well-understood and can be accurately approximated by linear models, especially models based on orbit elements \cite{Sullivan2017}.
%

Solution methodologies in literature for this class of problem can be divided into three broad categories: closed-form solutions, direct optimization methods, and indirect optimization methods.
Closed-form solutions are highly desirable because they are robust, predictable, and computationally efficient.
However, such solutions are inherently specific to the prescribed state representation, dynamics model, and cost function.
Indeed, despite decades of research, such solutions have only been found to date for a limited class of problems in spacecraft formation-flying \cite{DAmicoThesis,Gaias2015diffdrag1,Chernick2018,Chernick2018reachable,Serra2018}.

Direct optimization methods offer a greater degree of generality by formulating the optimal control problem as a nonlinear program with the times, magnitudes, and directions of the applied control inputs as variables \cite{Betts1998}.
The simplest approach to direct optimization is to discretize the admissible control window and optimize the applied control inputs at each of these times.
However, this approach requires enormous computational resources for all but the simplest problems.
To reduce computation effort, previous authors have developed iterative approaches that refine a small set of candidate control input times.
However, the minimum cost to reach a specified target state is a non-convex function of the number of impulses and the times at which they are applied \cite{Sobiesiak2015}.
As a result, these algorithms cannot guarantee convergence to a globally optimal solution \cite{Kim2002,Kim2009}.

Due to this weakness, the majority of numerical approaches in literature use indirect optimization techniques, which can be divided into two general approaches.
The first indirect optimization approach is based on some form of Lawden's so-called ``primer vector" \cite{Lawden1963}, which is an alias for the part of the costate that governs the control input according to Pontryagin's maximum principle.
Using this method, the optimal control problem is cast as a two-point boundary value problem where an optimal solution must satisfy a set of analytical conditions on the evolution of the primer vector.
This approach has been studied continuously for over fifty years \cite{Handelsman1968, Prussing1969, Prussing1970, Jezewski1980, Roscoe2015, Arzelier2016}, but the resulting algorithms are subject to substantial limitations.
For example, Roscoe's algorithm \cite{Roscoe2015} is known to have a limited radius of convergence from the initial estimate of the times of optimal control inputs.
Instead, the algorithm proposed by Arzelier \cite{Arzelier2016} provides guaranteed convergence to a globally optimal solution by sequentially adding candidate control input times based on the optimality criteria.
However, Arzelier's algorithm is developed under two limiting assumptions: 1) the cost of a control input is a constant $p$-norm function, and 2) the columns of the control input matrix are linearly independent.
Also, no considerations are made regarding the sensitivity of the cost of feasible solutions to errors in the control input times in corner cases.

The second widely studied indirect optimization approach is based on reachable set theory.
This approach was developed in the late 1960s with key contributions given by Neustadt \cite{Neustadt1964}, Barr \cite{Barr1969}, and Gilbert \cite{Gilbert1971}.
In contrast to algorithms based on primer vector theory, Gilbert developed an algorithm that provides global convergence for problems where the cost is expressed as a constant norm-like function of the control input \cite{Gilbert1971}.
%
%
Additionally, this algorithm is posed in a way that can be adapted to more challenging problems provided that the reachable set exhibits certain properties.
However, for some unknown reason extensions of this approach have not been studied in recent literature.

A common limitation of all of these algorithms is that the cost of a specified control input is not allowed to vary over time.
An algorithm that provides optimal solutions without this constraint could be applied to problems with complex, time-varying behaviors such a spacecraft with multiple attitude modes.
%

To meet this need, this paper makes contributions to the state-of-the-art in both theory and application.
The theoretical contribution includes derivation of the properties of the set of states that can be reached by admissible control input profiles with a specified cost.
Specifically, it is shown that the reachable sets for this problem and Gilbert's problem \cite{Gilbert1971} have identical properties.
With this in mind, contact and support functions for the reachable set are derived to enable use of existing algorithms for optimization over parameterized sets for this class of problem.
Next, the optimization problem is reformulated as a semi-infinite convex program and corresponding optimality conditions are derived.
It is shown that under the same additional assumptions, this semi-infinite convex program is identical to the form developed using primer vector theory \cite{Arzelier2016}, providing a geometric interpretation of the time-invariant portion of the primer vector.
The application contribution of this paper is a new three-step algorithm that quickly computes globally optimal impulsive control input sequences for the considered class of fuel-optimal control problems.
The three steps of this algorithm proceed as follows.
First, an initial set of candidate times for control inputs is computed from an a-priori estimate of the outward normal direction to the reachable set at the target state.
Second, the set of control input times and outward normal vector are iteratively refined until the optimality conditions are satisfied to within a user-specified tolerance.
To minimize computation cost, candidate times at which an optimal control cannot be applied are discarded in each iteration.
Additionally, each iteration provides a feasible solution with bounded sub-optimality, making the algorithm well-suited to real-time applications.
Third, an optimal sequence of impulsive control inputs is computed by solving a quadratic program.

The proposed algorithm is validated in four steps.
First, the performance of the algorithm is demonstrated through implementation in a challenging example problem based on the proposed Miniaturized Distributed Occulter/Telescope (mDOT) small satellite mission \cite{Koenig2019}.
Second, the computation time and number of required iterations for the proposed algorithm are compared with Gilbert's algorithm (using the contact and support functions developed in this paper) \cite{Gilbert1971} and direct optimization in Monte Carlo simulations.
Third, the Monte Carlo simulations with the proposed algorithm are repeated with two different initialization schemes to demonstrate that the algorithm is robust to poor initial conditions.
Finally, a selection of example problems are solved using widely varying discretizations of the time domain to characterize the sensitivity of the required computation time to the number of candidate control input times.
%

%
%
%
%


\section{Problem Definition}\label{sec:problemdefinition}

Consider a linear dynamical system with state vector $\boldsymbol x(t)\in\mathbb{R}^n$ and control input vector $\boldsymbol u(t)\in\mathbb{R}^{m}$ with dynamics that evolve as given by
\begin{equation}
\label{eq:diffdyn}
\dot{\boldsymbol x}(t) = \mathbf A(t)\boldsymbol x(t)+\mathbf B(t)\boldsymbol u(t)
\end{equation}
where $\mathbf A(t)\in\mathbb{R}^{n\times n}$ is the plant matrix and $\mathbf B(t)\in\mathbb{R}^{n\times m}$ is the control input matrix.
The only assumptions imposed on these matrices are that they are real and continuous on the closed interval [$t_i$, $t_f$], where $t_i$ denotes the initial time and $t_f$ denotes the final time.
Next, suppose that $\mathbf\Psi(t)$ is a fundamental matrix solution of the homogeneous equation ($\boldsymbol u(t) \equiv \mathbf 0 $) associated to (\ref{eq:diffdyn}).
Using this solution, a state transition matrix (STM) $\mathbf\Phi(t,t+\tau)$ that propagates the state from time $t$ to $t+\tau$ can be defined as
\begin{equation}
\label{eq:stmdef}
\mathbf\Phi(t,t+\tau) = \mathbf\Psi(t+\tau)\mathbf\Psi^{-1}(t).
\end{equation}
Using (\ref{eq:diffdyn}) and (\ref{eq:stmdef}), the final state $\boldsymbol x(t_f)$ can be expressed as a function of the initial state $\boldsymbol x(t_i)$ and the control input profile as given by
\begin{equation}
\label{eq:dynamics}
\begin{split}
\boldsymbol x(t_f) = \mathbf\Phi(t_i,t_f)\boldsymbol x(t_i)+\int_{t_i}^{t_f}\mathbf\Phi(\tau,t_f)\mathbf B(\tau)\boldsymbol u(\tau)d\tau.
\end{split}
\end{equation}
To simplify notation as in \cite{Carter1995}, let the pseudostate $\boldsymbol w$ and matrix $\mathbf\Gamma(t)$ be defined as
\begin{align*}
\boldsymbol w = \boldsymbol x(t_f)-\mathbf\Phi(t_i,t_f)\boldsymbol x(t_i), \quad \mathbf\Gamma(t) = \mathbf\Phi(t,t_f)\mathbf B(t).
\end{align*}

Using this notation, the class of optimal control problem considered in this paper is posed as follows
\begin{equation}
\label{eq:originalcontrolproblem}
\begin{split}
\textrm{minimize:} \quad J(\boldsymbol u(t)) = \int_{t_i}^{t_f} f(\boldsymbol u(\tau), \tau)d\tau \\
\textrm{subject to} \quad \boldsymbol w = \int_{t_i}^{t_f}\mathbf\Gamma(\tau)\boldsymbol u(\tau)d\tau \quad
\end{split}
\end{equation}
where the decision variable is the control input profile.
To simplify later discussions, it is hereafter assumed that $\boldsymbol w$ is reachable and nonzero.
This causes no loss in generality because the trivial solution $\boldsymbol u(t) = \boldsymbol 0$ is optimal for $\boldsymbol w = \boldsymbol 0$. 

The constraints for admissible control input profiles in this problem are as follows.
Nonzero control inputs are allowed at any time in the compact set $T \subseteq [t_i, t_f]$ and $T$ is divided into finite number $o$ of mutually exclusive compact sets $T_j$ such that
\begin{equation}
    T = \bigcup_{j = 1,\hdots,o} T_j, \qquad T_k \cap T_l = \varnothing \:\: \forall k \ne l.
\end{equation}
While this definition does require gaps between each $T_j$, this constraint does not compromise the generality of the solution approach for practical applications (i.e. this requirement is always satisfied if $T$ is a discrete set of times).
Within each interval $T_j$, the set of admissible control inputs is a closed convex cone $C_j \in \mathbb{R}^m$ with its vertex located at the origin.
This property allows for the possibility of inadmissible control inputs in $\mathbb{R}^m$ (e.g. directions in which a thruster cannot be fired).
Accordingly, the set of admissible control input profiles $\mathcal{U}$ is defined as
\begin{equation}
    \mathcal{U} = \begin{Bmatrix} \boldsymbol u(t): \boldsymbol u(t) \in \begin{Bmatrix} C_j & \forall t \in T_j,\; 1 \le j \le o \\ \boldsymbol 0 & \textrm{otherwise} \qquad \qquad \quad \end{Bmatrix} \end{Bmatrix}.
\end{equation}
The cost function $f(\boldsymbol u(\tau), \tau)$ is piecewise-defined over the intervals $T_j$ as given by
\begin{equation}
    \label{eq:piecewisecostfunction}
    f(\boldsymbol u(t), t) =
    \begin{Bmatrix}
    f_1(\boldsymbol u(t)), & t \in T_1 \\
    f_2(\boldsymbol u(t)), & t \in T_2 \\
    & \vdots & \\
    f_o(\boldsymbol u(t)), & t \in T_o \\
    0 & \textrm{otherwise}
    \end{Bmatrix}
\end{equation}
In this definition, each $f_j$ is a norm-like function with the three properties specified by Gilbert \cite{Gilbert1971} including: 1) $f_j$ is defined for all $\boldsymbol u \in C_j$, 2) the sublevel sets $U_j(c)$ defined as
\begin{equation}
    \label{eq:Udef}
    U_j(c) = 
    \begin{Bmatrix}
    \boldsymbol u : \boldsymbol u \in C_j, f_j(\boldsymbol u) \le c
    \end{Bmatrix}
\end{equation}
are convex and compact for all $c \ge 0$, and 3) $f_j(\alpha \boldsymbol u) = \alpha f_j(\boldsymbol u)$ for all $\alpha \ge 0$.
The last property ensures that the cost of a control input applied at a specified time scales linearly with its magnitude and that all nonzero control inputs have nonzero cost.
Under this definition, the set of admissible control inputs for Gilbert's problem \cite{Gilbert1971} is a subset of admissible control inputs in this problem where $o = 1$.

Some example norm-like cost functions are included in Table \ref{tab:costfunctions} along with operational constraints that provide the specified cost behavior.
Many of these cost functions are results of attitude constraints imposed on a spacecraft, which can vary over time, motivating the need for a solver that can accommodate a piecewise-defined cost function.
In the last example, the cost is the solution to the linear program given by
\begin{equation}
    \textrm{minimize:} \; \boldsymbol 1^T\boldsymbol\alpha \quad \textrm{subject to:} \; \boldsymbol u = \mathbf V^{vertex}\boldsymbol\alpha, \; \boldsymbol\alpha \ge 0.
\end{equation}
However, the sublevel sets of this cost function (the sets of $\boldsymbol u$ for which there exists a feasible $\boldsymbol\alpha$ that satisfies $\boldsymbol 1^T\boldsymbol\alpha \le c$ for some $c \ge 0$) are compact polyhedra with vertices parallel to columns of $\mathbf V^{vertex}$.
Thus, it is possible to compute a matrix $\mathbf V^{face} \in \mathbb{R}^{M\times m}$ for any $\mathbf V^{vertex}$ such that the sublevel set for cost $c$ is the set of solutions to $\mathbf V^{face}\boldsymbol u \le c\boldsymbol 1$.
Accordingly, the cost can be expressed as $\max(\mathbf V^{face}\boldsymbol u)$ where the function $\max$ returns the maximum of all elements in the vector argument.
Using this matrix, evaluation of the cost function is simplified to computing the maximum value of a vector.
\begin{table}[htb]
\caption{Example norm-like cost functions in $\mathbb{R}^3$.}
\centering
\begin{tabular}{m{2.5cm} m{5.5cm}}
\hline
Cost function & Operational constraints\\ \hline
$||\boldsymbol u||_2$ & Spacecraft that can align a single thruster with the desired maneuver direction. \\ \hline
$||\boldsymbol u||_1$ & Spacecraft with fixed attitude and three pairs of thrusters mounted on opposite sides on mutually perpendicular axes. \\ \hline
$|u_1|+\sqrt{u_2^2+u_3^2}$ & Spacecraft with two pairs of thrusters on perpendicular axes where one axis is fixed. \\ \hline
$\max(\mathbf V^{face}\boldsymbol u)$ & 
Spacecraft with fixed attitude and $N$ thrusters aligned with columns of $\mathbf V^{vertex} \in \mathbb{R}^{m\times N}$. Note: only valid for convex cone of the form $\boldsymbol u = \mathbf V^{vertex} \boldsymbol\alpha$ where $\boldsymbol\alpha$ is in the nonnegative orthant. \\  \hline
\end{tabular}
\label{tab:costfunctions}
\end{table}

Finally, Neustadt demonstrated that there must exist optimal control input profiles consisting of $n$ or fewer impulses for a subset of this class of problem where the cost function is a constant $p$-norm \cite{Neustadt1964}.
With this in mind, attention is hereafter restricted to impulsive control input profiles of the form
\begin{equation}
    \boldsymbol u(t) = \sum_{j=1}^N \delta(t-t_j)\boldsymbol v_j
\end{equation}
where $\delta$ denotes the Dirac delta function, $N$ is the number of applied impulses, $t_j$ are the times when each impulse is applied, and $\boldsymbol v_j$ defines the magnitude and direction of each impulse.
To avoid any misinterpretation of the meaning of the impulse function, the evolution of the state including an impulsive control input profile is given as
\begin{equation}
    \boldsymbol x(t) = \mathbf\Phi(t_i,t)\boldsymbol x(t_i) +\sum_{j=1}^{N(t)} \mathbf\Phi(t_j,t)\mathbf B(t_j)\boldsymbol v_j
\end{equation}
where $N(t)$ is the number of impulses such that $t_j \le t$.
The optimality of impulsive control input profiles for this problem class is proven in Section \ref{sec:reformulation}.

\section{Set Definitions and Properties}\label{sec:reachablesets}

The algorithm proposed in this paper is based on the properties of the set of pseudostates that can be reached by admissible control input profiles at a specified cost, which are derived in the following.
The following notation conventions are adopted throughout this derivation for simplicity.
First, the convex hull of $S$ is denoted co$S$.
Second, $\alpha S$ denotes the set defined as
\begin{equation}
    \alpha S = \begin{Bmatrix} \bar{\boldsymbol y} : \bar{\boldsymbol y} = \alpha\boldsymbol y, \; \boldsymbol y \in S \end{Bmatrix}
\end{equation}

It is instructive to first consider the set of pseudostates that can be reached with control inputs applied in only one interval $T_j$.
Let $\mathcal{U}^1_j(c)$ be the set of single admissible impulsive control inputs executed in $T_j$ with cost no greater than $c$, which is defined as
\begin{equation}
\label{eq:U1jcdef}
    \mathcal{U}^1_j(c) = \begin{Bmatrix} \boldsymbol u(t): & \boldsymbol u(t) = \delta(t-t_{cont})\boldsymbol v, \\ & t_{cont} \in T_j, \: \boldsymbol v \in U_j(c) \end{Bmatrix}.
\end{equation}
The set of pseudostates that can be reached by these impulses, denoted $S^1_j(c)$, is defined as 
\begin{equation}
\label{eq:S1jcdef}
    S^1_j(c) = \begin{Bmatrix} \boldsymbol y: \boldsymbol y = \int_{t_i}^{t_f}\mathbf\Gamma(\tau)\boldsymbol u(\tau)d\tau, \; \boldsymbol u(t) \in \mathcal{U}^1_j(c) \end{Bmatrix}.
\end{equation}
Similarly, let $\mathcal{U}_j(c)$ denote the set of admissible control input profiles consisting of a finite number $N$ of impulses applied in $T_j$ with a total cost no greater than $c$.
Using the linearity and convexity of the cost function, $\mathcal{U}_j(c)$ can be defined as
\begin{equation}
\label{eq:Ujcdef}
    \mathcal{U}_j(c) = \begin{Bmatrix} \boldsymbol u(t): & \boldsymbol u(t) = \displaystyle\sum_{j=1}^N \alpha_j\boldsymbol u_j(t), \: \alpha_j \ge 0, \\ & \displaystyle\sum_{j=1}^N\alpha_j = 1, \: \boldsymbol u_j(t) \in \mathcal{U}^1_j(c) \end{Bmatrix},
\end{equation}
which is simply the set of convex combinations of elements of $\mathcal{U}^1_j(c)$.
Accordingly, the set of pseudostates that can be reached by a set of impulses applied in $T_j$ with with a cost no greater than $c$, denoted $S_j(c)$, can be defined as
\begin{equation}
    \label{eq:Sjcdef}
    S_j(c) = \textrm{co} S^1_j(c).
\end{equation}
The set $S_j(c)$ is equivalent to the reachable set considered by Gilbert and its relationship with $S^1_j(c)$ is proven in Theorem 1 in \cite{Gilbert1971}.

With this relationship in mind, the reachable set will now be generalized to include all $T_j$.
%
%
%
The set of admissible single impulses with cost no greater than $c$ is denoted $\mathcal{U}^1(c)$ and is defined as
\begin{equation}
\label{eq:U1cdef}
    \mathcal{U}^1(c) = \begin{Bmatrix} \boldsymbol u(t): \; \boldsymbol u(t) \in \mathcal{U}^1_j(c), \: 1 \le j \le o \end{Bmatrix}
\end{equation}
and the set of pseudostates that can be reached by these impulses, denoted $S^1(c)$, is defined as
\begin{equation}
    \label{eq:S1cdef}
    S^1(c) = \begin{Bmatrix} \boldsymbol y:\; \boldsymbol y = \int_{t_i}^{t_f}\mathbf\Gamma(\tau)\boldsymbol u(\tau)d\tau, \; \boldsymbol u(t) \in \mathcal{U}^1(c) \end{Bmatrix}.
\end{equation}
Combining (\ref{eq:S1jcdef}), (\ref{eq:U1cdef}), and (\ref{eq:S1cdef}), $S^1(c)$ can be expressed as 
\begin{equation}
\label{eq:U1S1Union}
    S^1(c) = \bigcup_{j = 1,\hdots,o} S^1_j(c).
\end{equation}
Next, using the linearity and convexity of norm-like cost functions, the set of admissible impulsive control input profiles of cost $c$ or less, denoted $\mathcal{U}(c)$, is defined as
\begin{equation}
\label{eq:Ucdef}
\mathcal{U}(c) = \begin{Bmatrix} \boldsymbol u(t): & \boldsymbol u(t) = \displaystyle\sum_{j=1}^N \alpha_j\boldsymbol u_j(t), \: \alpha_j \ge 0, \\ & \displaystyle\sum_{j=1}^N\alpha_j = 1, \: \boldsymbol u_j(t) \in \mathcal{U}^1(c) \end{Bmatrix},
\end{equation}
which is the set of convex combinations of elements of $\mathcal{U}^1(c)$.
Accordingly, the set of pseudostates that can be reached by such a control input profile, denoted $S(c)$, is given by
\begin{equation}
\label{eq:Scdef2}
\begin{split}
    S(c) = \textrm{co}S^1(c).
\end{split}
\end{equation}
Moreover, by replacing $S_j(c)$ and $S^1_j(c)$ with $S(c)$ and $S^1(c)$ in Theorem 1 in \cite{Gilbert1971}, it is proven that $S(c)$ has the following properties:\\

\noindent\emph{Property 1:} $S(c)$ contains the origin for all $c \ge 0$.

\noindent\emph{Property 2:} $S(c_1) \subseteq S(c_2)$ for any $0 \le c_1 \le c_2$.

\noindent\emph{Property 3:} $S(c) = cS(1)$ for all $c \ge 0$.

\noindent\emph{Property 4:} $S(c)$ is convex, compact and Lipschitz continuous in the Hausdorff metric for all $c \ge 0$.

\noindent\emph{Property 5:} For any $\boldsymbol w \in S(c)$, there exists a control input profile that reaches $\boldsymbol w$ and consists of $n$ or fewer impulses with a total cost of $c$ or less.\\

This derivation has shown that for problems where the time domain can be expressed as a union of a finite number of compact sets, using a piecewise-defined norm-like cost function instead of a constant norm-like function has no impact on the properties of the reachable sets.
Accordingly, it is possible to use old algorithms in optimization over parameterized sets (e.g. Gilbert's algorithm \cite{Gilbert1971}) to solve this class of problem provided that the contact and support functions of $S(c)$ can be evaluated.


\section{Contact and Support Functions} \label{sec:contactsupport}

\begin{figure}[htb!]
\begin{center}
\includegraphics[trim={9cm 5cm 13cm 3.5cm},clip,width=0.6\columnwidth]{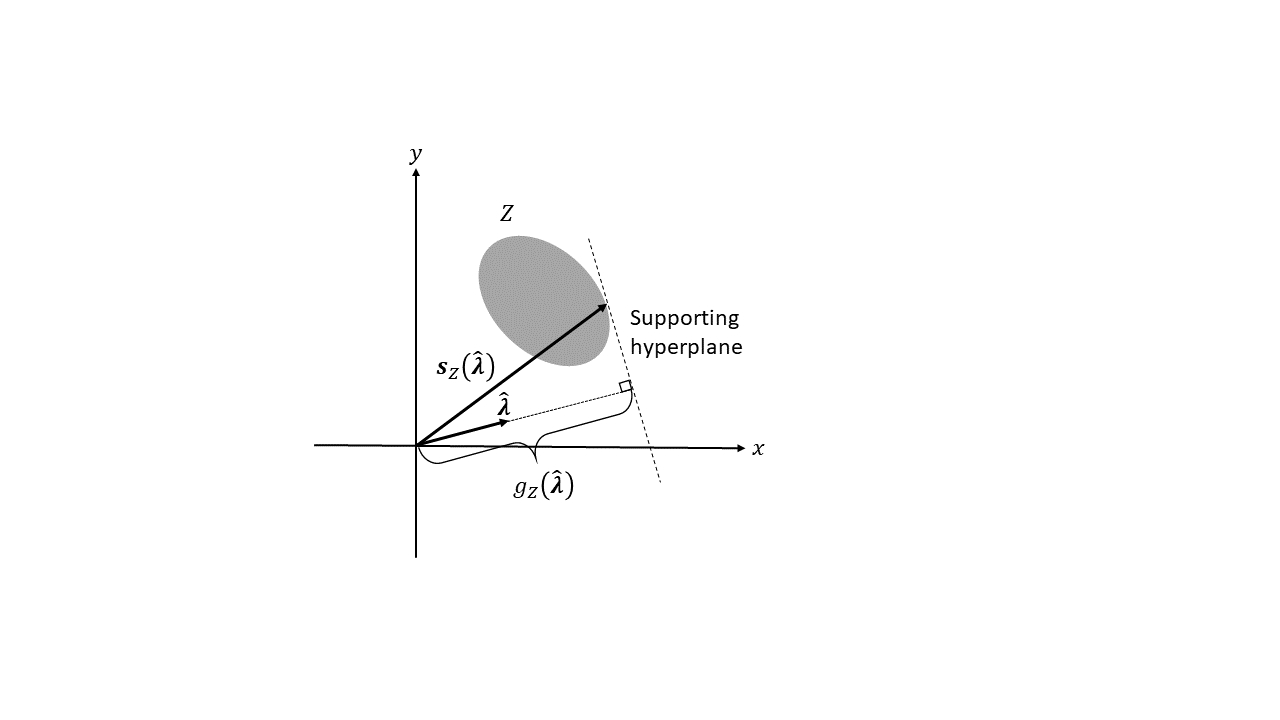}
\caption{Illustration of geometric relationships between $Z$ (gray) and its contact and support functions for a unit vector $\hat{\boldsymbol\lambda}$ and corresponding supporting hyperplane (dashed line).}\label{fig:contactsupport}
\end{center}
\end{figure}

To facilitate use of the algorithm proposed in Section \ref{sec:algorithm} and other existing algorithms in literature, contact and support functions for $S(c)$ are derived in the following.
A contact function for a set $Z \in \mathbb{R}^n$ is any function from $\mathbb{R}^n$ to $\mathbb{R}^1$ that satisfies
\begin{equation}
    \label{eq:contact}
    g_Z(\boldsymbol\lambda) = \Big[ \max_{\boldsymbol z \in Z} \boldsymbol\lambda^T\boldsymbol z \Big].
\end{equation}
where $\boldsymbol\lambda$ is any vector in $\mathbb{R}^n$.
Similarly, a support function $\boldsymbol s_Z(\boldsymbol\lambda)$ is a function from $\mathbb{R}^n$ to $Z$ such that
\begin{equation}
    \label{eq:support}
    \boldsymbol\lambda^T \boldsymbol s_Z(\boldsymbol\lambda) = g_Z(\boldsymbol\lambda).
\end{equation}
The geometric relationships between $Z$ and its contact and support functions for a unit vector $\hat{\boldsymbol\lambda}$ are illustrated in Fig. \ref{fig:contactsupport}.

The algorithm proposed in this paper and Gilbert's algorithm both require evaluation of the contact and support functions of $S(c)$.
Using the properties of the sets defined in the previous section, a simple procedure for these evaluations is developed in the following.
Using Property 3 of $S(c)$, the contact function can be expressed as
\begin{equation}
    \label{eq:sub_i}
    g_{S(c)}(\boldsymbol\lambda) = cg_{S(1)}(\boldsymbol\lambda) 
\end{equation}
Because the contact function of any closed set is equal to the contact function of its convex hull (Theorem 1 in \cite{Barr1969}), $g_{S(1)}$ can be simplified as
\begin{equation}
    g_{S(1)}(\boldsymbol\lambda) = g_{S^1(1)}(\boldsymbol\lambda).
\end{equation}
Because $S^1(j)$ is the union of a finite number of $S_j(1)$ (from (\ref{eq:U1S1Union})), $g_{S^1(1)}$ can be expressed as
\begin{equation}
   g_{S^1(1)}(\boldsymbol\lambda) = \max_{j = 1,\hdots,o} g_{S^1_j(1)}(\boldsymbol\lambda)
\end{equation}
Finally, using the definition of $S^1_j(1)$ in (\ref{eq:S1jcdef}), each $g_{S^1_j(1)}(\boldsymbol\lambda)$ can be formulated as
\begin{equation}
    \label{eq:sub_f}
   g_{S^1_j(1)}(\boldsymbol\lambda) =  \max_{t\in T_j} g_{U_j(1)}(\mathbf\Gamma^T(t)\boldsymbol\lambda)
\end{equation}
Combining (\ref{eq:contact}) and (\ref{eq:sub_i})-(\ref{eq:sub_f}), a contact function for $S(c)$ can be formulated as
\begin{equation}
    \label{eq:simplecontact}
    g_{S(c)}(\boldsymbol\lambda) =
    c\max_{j = 1,\hdots,o} \Big( \max_{t\in T_j} g_{U_j(1)}(\mathbf\Gamma^T(t)\boldsymbol\lambda)  \Big)
\end{equation}
To simplify this expression, let $U(1,t)$ be the sublevel set of unit cost at time $t$, which is defined as
\begin{equation}
    U(1,t) = \begin{Bmatrix} U_j(1) & t \in T_j, & 1 \le j \le o \\ \mathbf 0 & \textrm{otherwise} \end{Bmatrix}.
\end{equation}
Using this definition, (\ref{eq:simplecontact}) can be simplified to
\begin{equation}
    \label{eq:simplestcontact}
    g_{S(c)}(\boldsymbol\lambda) =
    c \max_{t\in T}  g_{U(1, t)}(\mathbf\Gamma^T(t)\boldsymbol\lambda) 
\end{equation}

Next, let $T^*$ be  the set of times at which $g_{U(1, t)}(\mathbf\Gamma^T(t)\boldsymbol\lambda)$ takes on its maximum value in the domain $T$, which is defined as
\begin{equation}
    T^* = \begin{Bmatrix} t: & g_{U(1, t)}(\mathbf\Gamma^T(t)\boldsymbol\lambda) = g_{S(1)}(\boldsymbol\lambda), & t \in T \end{Bmatrix}.
\end{equation}
Using this definition, a support function $\boldsymbol s_{S(c)}(\boldsymbol\lambda)$ can be formulated as
\begin{equation}
    \label{eq:simplesupport}
    \boldsymbol s_{S(c)}(\boldsymbol\lambda) = c\mathbf\Gamma(t^*) \boldsymbol s_{U(1,t^*)}(\mathbf\Gamma^T(t^*)\boldsymbol\lambda), \quad t^* \in T^*.
\end{equation}

It is evident from (\ref{eq:simplestcontact}) and (\ref{eq:simplesupport}) that evaluating $g_{S(c)}(\boldsymbol\lambda)$ and $\boldsymbol s_{S(c)}(\boldsymbol\lambda)$ requires a global optimization over the time domain $T$.
If $T$ is a set of discrete samples, the computation cost varies linearly with the number of samples.
It is also evident that the computation cost of these evaluations may be impractically large if contact and support functions to $U(1,t)$ cannot be evaluated analytically.
However, it is expected that analytical contact and support functions can be derived for a wide range of cost functions of interest.
For example, Table \ref{tab:examplecontactsupport} includes contact and support functions for all the example cost functions in Table \ref{tab:costfunctions}.
The contact functions for three of these examples are evaluated by computing the largest inner product between $\boldsymbol\lambda$ and an element in $\mathbf W_j$, which is the set of columns of the matrix $\mathbf W$ (shown in the second column of the table when relevant).
The corresponding support functions provide the column(s) of $\mathbf W$ that provides this maximum inner product.
\begin{table*}[htb]
\caption{Analytical support and contact functions for sublevel sets of example cost functions.}
\centering
\begin{tabular}{@{}cccc@{}}
\hline
$f_j(\boldsymbol u)$ & $\mathbf W$ & $\boldsymbol s_{U_j(1})(\boldsymbol\lambda)$ & $g_{U_j(1)}(\boldsymbol\lambda)$ \\ \hline
$||\boldsymbol u||_2$ & - & $\boldsymbol\lambda/||\boldsymbol\lambda||_2$ & $||\boldsymbol\lambda||_2$ \\ \hline
$||\boldsymbol u||_1$ & $\begin{bmatrix} \mathbf I & -\mathbf I \end{bmatrix}$ & $\displaystyle\argmax_{\boldsymbol w \in \mathbf W_j} \boldsymbol\lambda^T\boldsymbol w$ & $||\boldsymbol\lambda||_\infty$ \\ \hline
$|u_1|+\sqrt{u_2^2+u_3^2}$ & $\begin{bmatrix} \textrm{sign}(\lambda_1) & 0 \\ 0 & \frac{\lambda_2}{\sqrt{\lambda_2^2 + \lambda_3^2}} \\ 0 & \frac{\lambda_3}{\sqrt{\lambda_2^2 + \lambda_3^2}} \end{bmatrix}$ & $\displaystyle\argmax_{\boldsymbol w \in \mathbf W_j} \boldsymbol\lambda^T\boldsymbol w$ & $\displaystyle\max_{\boldsymbol w \in \mathbf W_j} \boldsymbol\lambda^T\boldsymbol w$ \\ \hline
$\max(\mathbf V^{face}\boldsymbol u)$ & $\begin{bmatrix} \mathbf 0^{m\times 1} & \mathbf V^{vertex} \end{bmatrix}$ & $\displaystyle\argmax_{\boldsymbol w \in \mathbf W_j} \boldsymbol\lambda^T\boldsymbol w$ & $\displaystyle\max_{\boldsymbol w \in \mathbf W_j} \boldsymbol\lambda^T\boldsymbol w$ \\  \hline
\end{tabular}\label{tab:examplecontactsupport}
\end{table*}


\section{Reformulation of the Optimization Problem}\label{sec:reformulation}

Using the properties of the sets defined in Section \ref{sec:reachablesets} and the contact functions derived in Section \ref{sec:contactsupport}, the optimal control problem in (\ref{eq:originalcontrolproblem}) is reformulated into a more computationally tractable form in the following.
Since the target pseudostate must be in the reachable set, an equivalent problem to (\ref{eq:originalcontrolproblem}) can be posed as
\begin{equation}
\label{eq:setdefprimal}
\textrm{minimize:} \quad c \qquad
\textrm{subject to:} \quad \boldsymbol w \in S(c).
\end{equation}
where the decision variable is the cost $c$.
Since $S(c)$ is compact and convex, it can be expressed as the intersection of all half-spaces that contain it (see Corollary 1 in Chapter 2 in \cite{Bazaraa1993}).
Thus, (\ref{eq:setdefprimal}) can be reformulated as
\begin{equation}
\label{eq:primalprob}
\begin{split}
\textrm{minimize:} \: c \qquad \qquad \qquad  \\
\textrm{subject to:} \: 
g_{S(c)}(\boldsymbol\lambda) \ge \boldsymbol\lambda^T\boldsymbol w \quad \forall \boldsymbol\lambda \in \mathbb{R}^n 
\end{split}
\end{equation}
%
%
%
Since $\boldsymbol w$ is reachable and nonzero, the constraint in (\ref{eq:primalprob}) is satisfied for any $\boldsymbol\lambda^T\boldsymbol w \le 0$ because $S(c)$ contains the origin for $c \ge 0$.
With this in mind, let $\boldsymbol\Lambda^+$ be defined as
\begin{equation}
    \boldsymbol\Lambda^+ = \begin{Bmatrix} \boldsymbol\lambda: & \boldsymbol\lambda^T\boldsymbol w > 0 \end{Bmatrix}
\end{equation}
Replacing $\forall\boldsymbol\lambda \in \mathbb{R}^n$ with $\forall\boldsymbol\lambda \in \boldsymbol\Lambda^+$ in (\ref{eq:primalprob}) results in no loss of generality and ensures that $g_{S(c)}(\boldsymbol\lambda) > 0$ for all admissible $\boldsymbol\lambda$.
Using this substitution and the linearity of the contact function with $c$, (\ref{eq:primalprob}) can be expressed as
\begin{equation}
\label{eq:primalprob2}
\begin{split}
\textrm{minimize:} \: c \qquad \qquad \quad  \\
\textrm{subject to:} \: 
c \ge \frac{\boldsymbol\lambda^T\boldsymbol w}{g_{S(1)}(\boldsymbol\lambda)} \qquad \forall \boldsymbol\lambda \in \boldsymbol\Lambda^+
\end{split}
\end{equation}
%
%
By inspection, the optimal solution $c^*$ to (\ref{eq:primalprob2}) is given by
\begin{equation}
    \label{eq:cstardef}
    c^* = \Big [  \max_{\boldsymbol\lambda \in \boldsymbol\Lambda^+} \frac{\boldsymbol\lambda^T\boldsymbol w}{g_{S(1)}(\boldsymbol\lambda)} \Big]
\end{equation}
and the set of feasible solutions to (\ref{eq:primalprob2}) is given by $c \ge c^*$.
From Theorem 8 in \cite{Barr1969}, there must exist a non-empty set $\boldsymbol\Lambda^* \subset \boldsymbol\Lambda^+$ such that 
\begin{equation}
\label{eq:Lambdastardef}
\begin{split}
    \frac{\boldsymbol\lambda^{*T}\boldsymbol w}{g_{S(1)}(\boldsymbol\lambda^*)} = \max_{\boldsymbol\lambda \in \boldsymbol\Lambda^+} \frac{\boldsymbol\lambda^T\boldsymbol w}{g_{S(1)}(\boldsymbol\lambda)} \quad \forall \boldsymbol\lambda^* \in \boldsymbol\Lambda^*\\
    \frac{\boldsymbol\lambda'^{T}\boldsymbol w}{g_{S(1)}(\boldsymbol\lambda')} < \max_{\boldsymbol\lambda \in \boldsymbol\Lambda^+} \frac{\boldsymbol\lambda^T\boldsymbol w}{g_{S(1)}(\boldsymbol\lambda)} \quad \forall \boldsymbol\lambda' \notin \boldsymbol\Lambda^*
\end{split}
\end{equation}
It is clear from this definition that $\boldsymbol\Lambda^*$ is the set of outward normal directions to supporting hyperplanes to $S(c)$ that contain $\boldsymbol w$.
With this in mind, consider the problem posed as
\begin{equation}
\label{eq:dualprob}
\begin{split}
\textrm{maximize:} \: c \qquad 
\textrm{subject to:} \: 
c \le \frac{\boldsymbol\lambda^T\boldsymbol w}{g_{S(1)}(\boldsymbol\lambda)} \qquad \boldsymbol\lambda \in \boldsymbol\Lambda^+
\end{split}
\end{equation}
where the decision variables are $c$ and $\boldsymbol{\lambda}$.\\

\noindent\emph{Theorem 1:}  If $c^{opt}$ and $\boldsymbol\lambda^{opt}$ are an optimal solution to (\ref{eq:dualprob}), then $c^{opt}$ is the optimal solution to (\ref{eq:primalprob}).\\

\noindent\emph{Proof:} By inspection of (\ref{eq:dualprob}) it is evident that $c^{opt}$ is equal to $c^*$ as defined in (\ref{eq:cstardef}), which means $c^{opt}$ is a feasible solution to (\ref{eq:primalprob}).
%
%
Any $c < c^{opt}$ cannot be a feasible solution to (\ref{eq:primalprob})  because it would violate the first constraint evaluated at $\boldsymbol\lambda^{opt}$.
Thus, $c^{opt}$ must be the optimal solution to (\ref{eq:primalprob}) $\blacksquare$.\\

\noindent\emph{Remark 1:} For any $\boldsymbol\lambda \in \boldsymbol\Lambda^+$, $\boldsymbol\lambda^T\boldsymbol w/g_{S(1)(\boldsymbol\lambda)}$ is a lower bound on the optimal objective in (\ref{eq:primalprob2}).
This property is leveraged in  the algorithm in Section \ref{sec:algorithm} to ensure that the converged solution has a cost within a specified tolerance of the optimal value.

Using the substitution in (\ref{eq:simplestcontact}) and simple arithmetic manipulations, (\ref{eq:dualprob}) can be reformulated as
\begin{equation}
\label{eq:dualprob2}
\begin{split}
\textrm{maximize:} \: \boldsymbol\lambda^T\boldsymbol w \qquad \qquad \\
\textrm{subject to:} \:  \max_{t\in T} g_{U(1,t)}(\mathbf\Gamma^T(t)\boldsymbol\lambda) \le 1
\end{split}
\end{equation}
where the decision variable is $\boldsymbol\lambda$ and the cost is encoded as $\boldsymbol\lambda^T\boldsymbol w$.
While this form is similar to the semi-infinite convex programs considered in numerous prior works (e.g. \cite{Arzelier2016} and \cite{Neustadt1964}), it is distinct in that it holds for an arbitrary piecewise-defined norm-like cost function instead of a constant $p$-norm cost.
This form of the optimal control problem is used in the solution algorithm proposed in Section \ref{sec:algorithm}.
%

From the properties of the optimal solutions to (\ref{eq:dualprob}) or (\ref{eq:dualprob2}) and the set definitions in Section \ref{sec:reachablesets}, necessary and sufficient optimality conditions for control input profiles are as follows.\\

\noindent\emph{Theorem 2:} Let $\boldsymbol\lambda^*$ be an optimal solution to (\ref{eq:dualprob2}) with corresponding optimal cost $c^* = \boldsymbol\lambda^{*T}\boldsymbol w$.
Additionally, let $T^*$ and $\mathcal{U}^*$ be defined as
\begin{equation}
\label{eq:optimalsetdefs}
\begin{split}
    T^* = \begin{Bmatrix} t: \: g_{U(1,t)}(\mathbf\Gamma^T(t)\boldsymbol\lambda^*) = 1,\: t \in T \end{Bmatrix} \quad \\
    \mathcal{U}^* = \begin{Bmatrix} \boldsymbol u(t): & \boldsymbol u(t) = \delta(t-t^*)\boldsymbol v, \: t^* \in T^* \\ & \boldsymbol v = c^*\boldsymbol s_{U(1,t^*)}(\mathbf\Gamma(t^*)\boldsymbol\lambda^*)\end{Bmatrix} 
\end{split}
\end{equation}
Any control input profile $\boldsymbol u^*(t)$ that satisfies
\begin{equation}
\label{eq:optimalcontrolprofile}
\begin{split}
    \boldsymbol u^*(t) = \sum_{j=1}^N \alpha_j\boldsymbol u_j(t), \; \boldsymbol u_j(t) \in \mathcal{U}^*, \; \alpha_j \ge 0, \\
    \sum_{j=1}^p \alpha_j = 1, \; \boldsymbol w = \int_{t_i}^{t_f}\mathbf\Gamma(t)\boldsymbol u^*(t) dt \qquad
\end{split}
\end{equation}
for positive integer $N$ is an optimal control input profile.\\

\begin{figure}[htb]
\begin{center}
\includegraphics[trim={6cm 5cm 6cm 5cm},clip,width=\columnwidth]{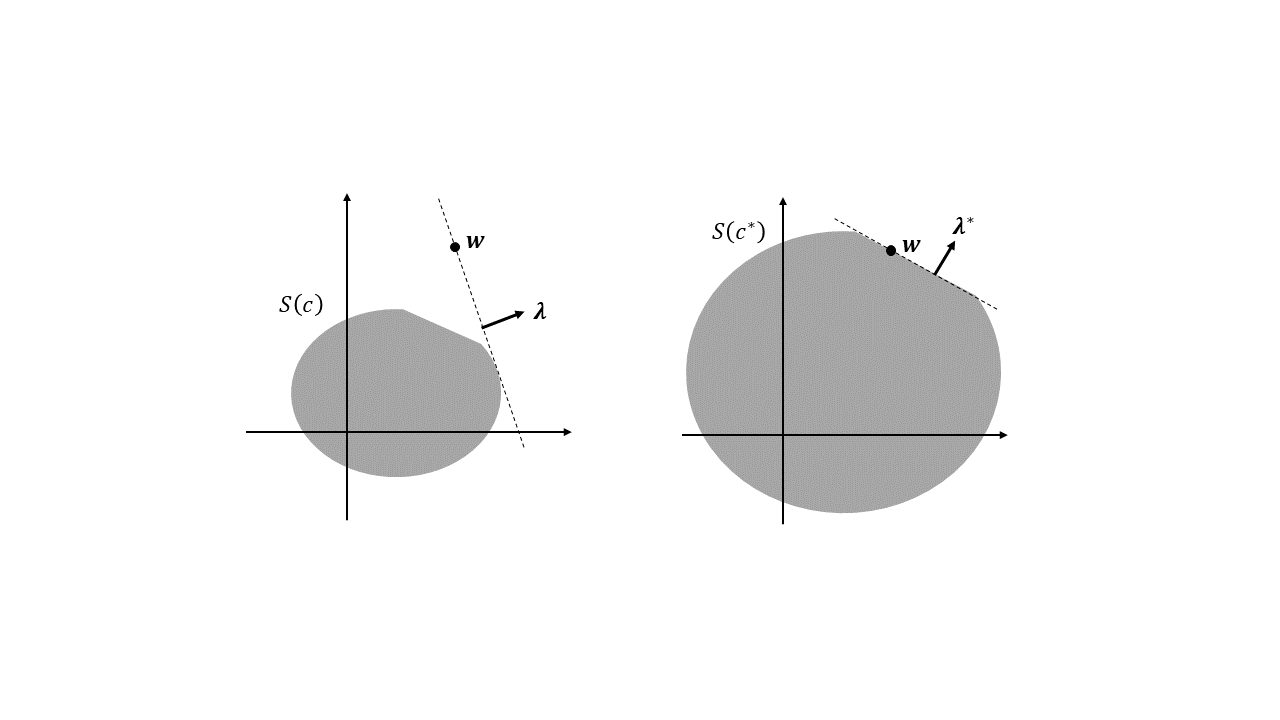}
\caption{Relationships between $\boldsymbol\lambda$, $\boldsymbol w$, and $S(c)$ for sub-optimal (left) and optimal (right) solutions to (\ref{eq:dualprob}) or (\ref{eq:dualprob2}).}\label{fig:lambdac}
\end{center}
\end{figure}

These properties are simple geometric consequences of the set definitions provided in Section \ref{sec:reachablesets} and can be obtained by replacing $\mathcal{U}^1_j(c^*)$ with $\mathcal{U}^1(c^*)$ in Theorem 2 in \cite{Gilbert1971}.
The geometric relationships between $\boldsymbol\lambda$, $\boldsymbol w$, and the corresponding $S(c)$ are illustrated in Fig. \ref{fig:lambdac} for sub-optimal (left) and optimal (right) solutions to (\ref{eq:dualprob}) or (\ref{eq:dualprob2}).

\noindent\emph{Remark 2:} Under the definitions in (\ref{eq:optimalsetdefs}), $T^*$ is the set of times at which it is possible to reach the hyperplane that is perpendicular to $\boldsymbol\lambda^*$ and contains $\boldsymbol w$ using a single impulse with a cost of $c^*$.
Similarly, $\mathcal{U}^*$ is the set of impulses of cost $c^*$ that reach this supporting hyperplane.

\noindent\emph{Remark 3:} From Property 5 of $S(c)$, there must exist at least one $\boldsymbol u^*(t)$ with $N \le n$.

\noindent\emph{Remark 4:} It is evident from (\ref{eq:optimalsetdefs}) and (\ref{eq:optimalcontrolprofile}) that an optimal control input profile can be computed from $\boldsymbol\lambda^*$ and $c^*$ using a simple three-step procedure.
First, $T^*$ is computed.
Second, for each time in $T^*$, candidate impulsive control inputs are computed by computing all values of $\boldsymbol s_{U(c^*,t)}(\boldsymbol\lambda^*)$.
Third, a convex combination of these impulses is computed that reaches $\boldsymbol w$.

Lastly, it has thus far been assumed that the control input profile consists of a set of impulses.
It will be proven in Theorem 3 that the cost associated with an optimal solution to (\ref{eq:dualprob2}) is optimal for all continuous control input profiles.\\

\noindent\emph{Theorem 3:} If $\boldsymbol\lambda^*$ is an optimal solution to (\ref{eq:dualprob2}) and $c^*$ is the corresponding cost, then $c^*$ is the minimum cost to reach $\boldsymbol w$ using any control input profile.\\

\noindent\emph{Proof:} Suppose there exists a control input profile $\boldsymbol u'(t)$ that reaches $\boldsymbol w$ at a cost of $c' < c^*$.
This implies that $\boldsymbol u'(t)$ satisfies
\begin{equation}
    \label{eq:hyperplaneconstraint}
    \boldsymbol\lambda^{*T}\boldsymbol w = c^* = \int_{t_i}^{t_f}\boldsymbol\lambda^{*T}\mathbf\Gamma(\tau)\boldsymbol u'(\tau)d\tau
\end{equation}
However, due to the linearity of the cost function, the integral must satisfy
\begin{equation}
    \int_{t_i}^{t_f}\boldsymbol\lambda^{*T}\mathbf\Gamma(\tau)\boldsymbol u'(\tau)d\tau \le c'\max_{t\in T}g_{U(1,t)}(\mathbf\Gamma(t)\boldsymbol\lambda^*)
\end{equation}
Using the constraint in (\ref{eq:dualprob2}), this simplifies to
\begin{equation}
    \int_{t_i}^{t_f}\boldsymbol\lambda^{*T}\mathbf\Gamma(\tau)\boldsymbol u'(\tau)d\tau \le c'
\end{equation}
which contradicts (\ref{eq:hyperplaneconstraint}).
Thus, there are no control input profiles with a cost less than $c^*$ that reach $\boldsymbol w$.
It follows that $c^*$ is the minimum cost to reach $\boldsymbol w$ using any control input profile $\blacksquare$.

\section{Comparison with Primer Vector Theory}\label{sec:primervector}

It is now worthwhile to consider the relationship between (\ref{eq:dualprob2}) and the formulations developed using primer vector theory \cite{Arzelier2016,Neustadt1964}.
These approaches assume that the cost function is a constant $p$-norm of the form
\begin{align*}
    f(\boldsymbol u) = ||\boldsymbol{u}||_p = \begin{Bmatrix} \Big( \sum_{j=1}^m |u_j|^{p} \Big)^{1/p}, & p \in [1, \infty) \\
    \displaystyle\max_{j = 1,\hdots,m}|u_j|, & p = \infty \end{Bmatrix}.
\end{align*}
Using H\"{o}lder's inequality, the contact function for the unit sublevel set of any $p$-norm can be expressed as
\begin{equation}
    g_{U(1,t)}(\mathbf\Gamma^T(t)\boldsymbol\lambda) = \max_{||\boldsymbol u||_p \le 1} \boldsymbol\lambda ^{T}\mathbf\Gamma(t)\boldsymbol u = ||\mathbf\Gamma^T(t)\boldsymbol\lambda||_q
\end{equation}
where $q$ is selected to satisfy $1/p + 1/q = 1$ and $\mathbf\Gamma^T(t)\boldsymbol\lambda$ is the primer vector.
Substituting this expression into (\ref{eq:dualprob2}) yields an equivalent form to Problem 5 in \cite{Arzelier2016}.
It is evident from this relationship that $g_{U(1,t)}(\mathbf\Gamma^T(t)\boldsymbol\lambda)$ is effectively a generalization of the $q$-norm of the primer vector for this class of problem.
However, unlike the magnitude of the primer vector, this function is not required to have a continuous first derivative with respect to time.
Additionally, it provides a geometric intuition that the time-invariant portion of the primer vector ($\boldsymbol\lambda$) is an outward normal direction to the reachable set at the location of the target pseudostate.
This property can be exploited by initializing solution algorithms with a reasonable a-priori estimate of the outward normal direction to the reachable set at the target state.


\section{A Fast and Robust Solution Algorithm}\label{sec:algorithm}

An efficient and robust algorithm that provides globally optimal impulsive control input sequences for the class of optimal control problems described in Section \ref{sec:problemdefinition} is presented in the following.
This algorithm includes three steps: 1) initialization, 2) iterative refinement, and 3) extraction of optimal control inputs.
The geometric properties of the problem described in the previous sections are leveraged at every step to minimize computation cost and maximize robustness to numerical errors.

\subsection*{Initialization}

The first step is generation of a set of candidate times $T^{est}$.  
The only requirement imposed on this step is that it must be possible to reach $\boldsymbol w$ using a combination of admissible impulses executed at times in $T^{est}$.
For most applications, a coarse discretization of $T$ is sufficient to meet this requirement.
However, to minimize the computation cost of refining the initial estimate, it is desirable to select candidates that are as close as possible to the optimal times for control inputs.

This can be accomplished by using an a-priori estimate of the optimal $\boldsymbol\lambda$, denoted $\boldsymbol\lambda_{est}$.
From the structure of the objective in (\ref{eq:dualprob}), a reasonable choice of $\boldsymbol\lambda_{est}$ is a vector parallel to $\boldsymbol w$.
Using this estimate, an initial set of candidate times for control inputs can be computed as follows.
First, a set of discrete samples $T^d \subseteq T$ is computed.
Next, $g_{U(1,t)}(\mathbf\Gamma^T(t)\boldsymbol\lambda_{est})$ is computed for each $t \in T^d$.
The initial set of control input times $T^{est}$ consists of the user-specified number $N$ of samples in $T^d$ for which $g_{U(1,t)}(\mathbf\Gamma^T(t)\boldsymbol\lambda_{est})$ is largest.
This initialization approach is summarized in the following pseudocode.\\

\noindent\textbf{Algorithm 1: Initialization}\\
\textbf{Inputs:} $T^d$, $\boldsymbol \lambda_{est}$, $\mathbf\Gamma(t)$, $N$\\
\textbf{Outputs:} $T^{est}$\\
\textbf{loop} $t \in T^d$\\
\indent compute $g_{U(1,t)}(\mathbf\Gamma^T(t)\boldsymbol\lambda_{est})$\\
\textbf{loop} $t \in T^d_j$\\
\indent \textbf{if} $g_{U(1,t)}(\mathbf\Gamma^T(t)\boldsymbol\lambda_{est})$ is one of the $N$ largest values\\
\indent \indent add $t$ to $T^{est}$ \\
\textbf{return} $T^{est}$\\

\begin{figure}[htb]
\begin{center}
\includegraphics[trim={11cm 5.5cm 11cm 5cm},clip,width=0.6\columnwidth]{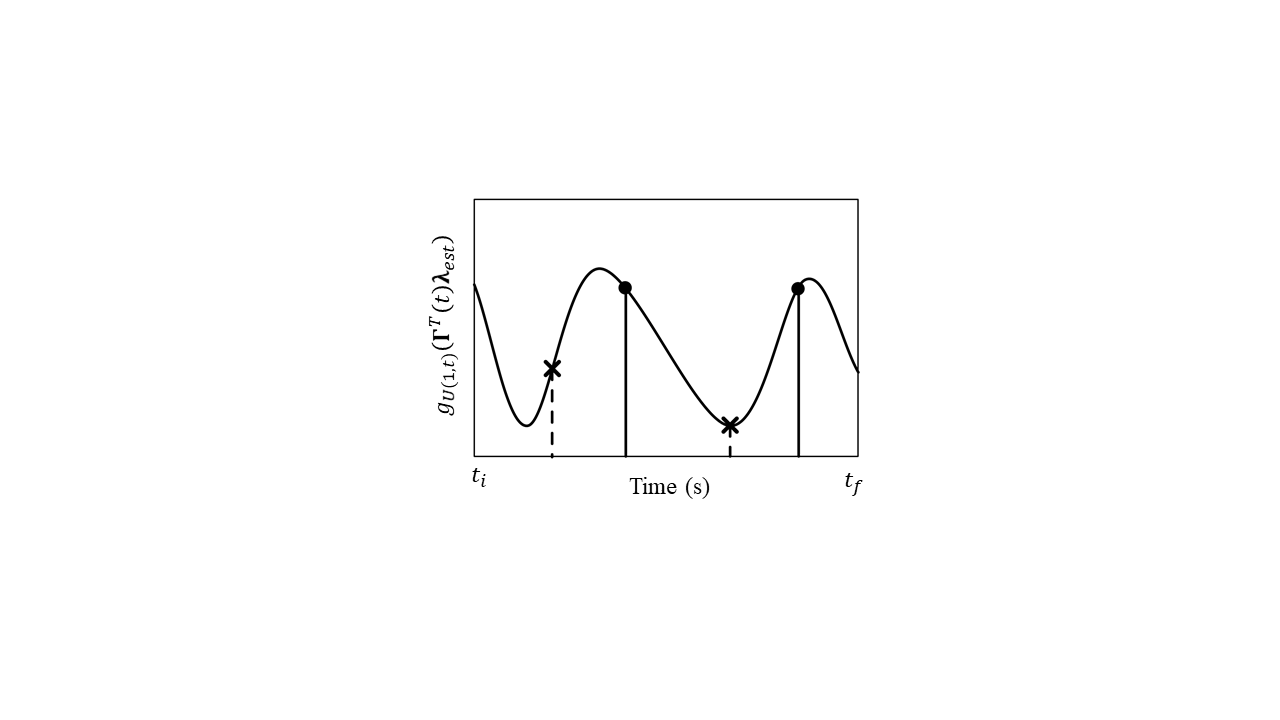}
\caption{Illustration of selection of two initial candidate times including selected times (circle) and rejected times (x).}\label{fig:initialization}
\end{center}
\end{figure}

A notional example of this initialization procedure is shown in Fig. \ref{fig:initialization}.
In this example, $T^d$ includes four times (indicated by vertical lines), and the initialization algorithm selects the two best times.
The two selected times (indicated by circles) are those at which $g_{U(1,t)}(\mathbf\Gamma^T(t)\boldsymbol\lambda_{est})$ is largest.
The rejected candidates are indicated by x markers.
%
%
%
%

\subsection*{Iterative Refinement}

Next, $\boldsymbol\lambda_{est}$ and $T^{est}$ are iteratively refined until they converge to $\boldsymbol\lambda_{opt}$ and $T^{opt}$, which satisfy the optimality criteria to within a user-specified tolerance.
This is accomplished using an iterative three-step procedure that provides global convergence from any feasible initial guess.
This procedure is inspired by Arzelier's algorithm \cite{Arzelier2016}, but includes modifications to simultaneously minimize the number of constraints that must be enforced in the required optimization problems and reduce the number of required iterations.
%
%
The first step in each iteration is computation of an optimal $\boldsymbol\lambda_{est}$ for the current instance of $T^{est}$.
This is accomplished by solving (\ref{eq:dualprob2}) with the modification that control inputs are only allowed at times in $T^{est}$.
This problem can be solved using conventional convex solvers, but the complexity of the problem depends on how the cost function is defined.
%

The second step is refinement of $T^{est}$ using the updated $\boldsymbol\lambda_{est}$.
This update procedure removes times at which optimal control inputs cannot be applied and adds candidate times that can reduce the total cost.
Specifically, all $t \in T^{est}$ that satisfy
\begin{align*}
g_{U(1,t)}(\mathbf\Gamma^T(t)\boldsymbol\lambda_{est}) < 1-\epsilon_{remove}
\end{align*}
for user-specified tolerance $\epsilon_{remove} > 0$ are removed from $T^{est}$ to reduce the number of constraints that must be enforced in subsequent iterations, thereby reducing computational effort.
Removing these times has no impact on the cost because optimal control inputs cannot be applied at these times.
Next, $g_{U(1,t)}(\mathbf\Gamma^T(t)\boldsymbol\lambda_{est})$ is evaluated for all times in $T$.
All times of local maxima of $g_{U(1,t)}(\mathbf\Gamma^T(t)\boldsymbol\lambda_{est})$ that are greater than one are added to $T^{est}$.
Adding these times ensures that $\max_{t\in T}g_{U(1,t)}(\mathbf\Gamma^T(t)\boldsymbol\lambda_{est})$ monotonically decreases to one with each subsequent iteration, thereby ensuring that $\boldsymbol\lambda_{est}$ and $T^{est}$ converge to $\boldsymbol\lambda_{opt}$ and $T^{opt}$, respectively.
While no rigorous guarantee is provided for the speed of convergence (which depends on the dynamics model and cost functions), the results in Section \ref{sec:validation} demonstrate that a wide range of problems can be solved in less than ten iterations.

The third step is to check the if the optimality criteria are satisfied to within a user-specified tolerance $\epsilon_{cost} > 0$.
Specifically, if the condition given by
\begin{align*}
\max_{t\in T}g_{U(1,t)}(\mathbf\Gamma^T(t)\boldsymbol\lambda_{est}) \le 1+\epsilon_{cost}
\end{align*}
is satisfied, then the algorithm terminates.
Otherwise, another iteration is performed.
This ensures that the cost of the final solution is within a factor of $\epsilon_{cost}$ of the optimal cost.
The described iteration procedure is summarized in the following pseudocode.\\

\noindent\textbf{Algorithm 2: Iterative Refinement}\\
\textbf{Inputs:} $T^{est}$, $\boldsymbol{\lambda}_{est}$, $T$, $\boldsymbol w$, $\mathbf\Gamma(t)$, $\epsilon_{cost}$, and $\epsilon_{remove}$\\
\textbf{Outputs:} $T^{opt}$ and $\boldsymbol\lambda_{opt}$\\
\textbf{do}\\
\indent $\boldsymbol\lambda_{est} \gets$ solution of problem:\\
\indent maximize: $\boldsymbol \lambda^T\boldsymbol w$\\
\indent subject to: $\displaystyle\max_{t\in T_{est}}g_{U(1,t)}(\mathbf\Gamma^T(t)\boldsymbol\lambda) \le 1$\\
\indent \textbf{loop} $t \in T_{est}$\\
\indent \indent \textbf{if} $g_{U(1,t)}(\mathbf\Gamma^T(t)\boldsymbol\lambda_{est}) < 1-\epsilon_{remove}$\\
\indent \indent \indent remove $t$ from $T^{est}$\\
\indent \textbf{loop} local maxima of $g_{U(1,t)}(\mathbf\Gamma^T(t)\boldsymbol\lambda_{est})$ in $T$\\
\indent \indent \textbf{if}  $g_{U(1,t)}(\mathbf\Gamma^T(t)\boldsymbol\lambda_{est}) > 1$\\
\indent \indent \indent add $t$ to $T^{est}$\\
\textbf{while}  $\displaystyle\max_{t\in T}g_{U(1,t)}(\mathbf\Gamma^T(t)\boldsymbol\lambda_{est}) > 1+\epsilon_{cost}$\\
$T_{opt} \gets T^{est}$\\
$\boldsymbol\lambda_{opt} \gets \boldsymbol\lambda_{est}$\\
\textbf{return} $T^{opt}$ and $\boldsymbol\lambda_{opt}$\\

\begin{figure}[htb]
\begin{center}
\includegraphics[trim={9cm 5cm 9cm 5cm},clip,width=0.8\columnwidth]{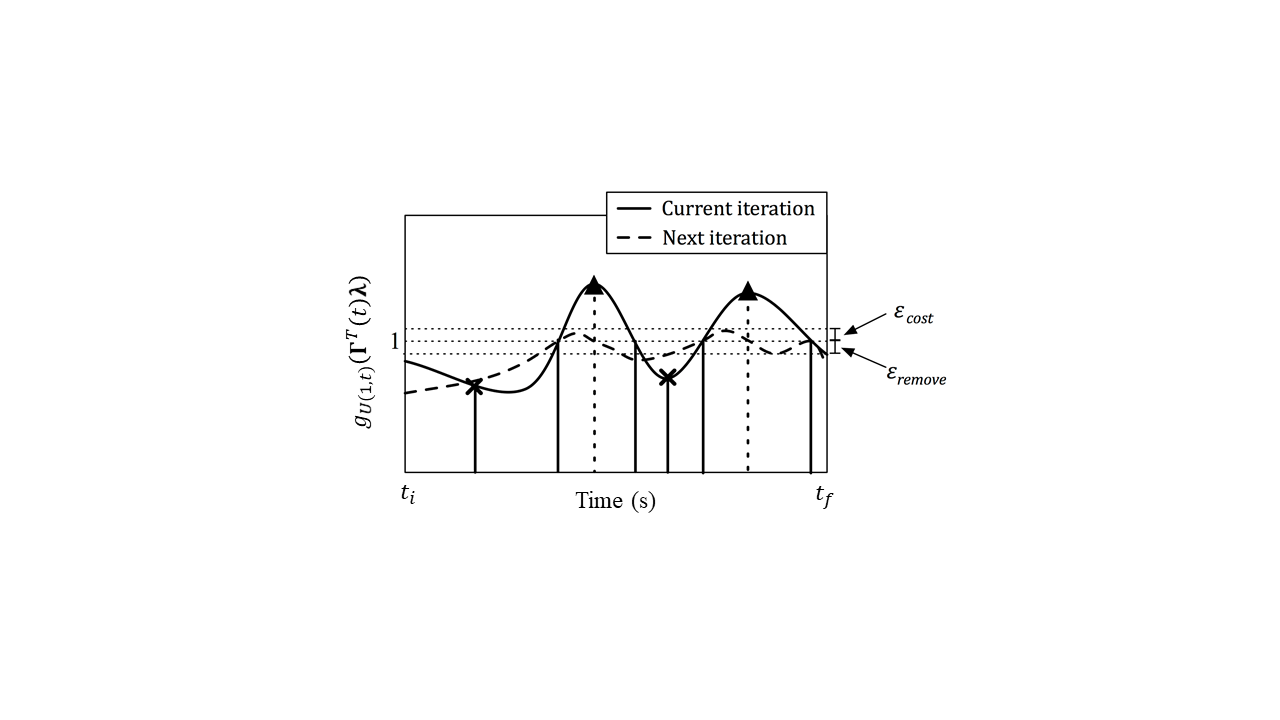}
\caption{Illustration of iterative refinement procedure including removed times (x) and added times (triangles).}\label{fig:refinement}
\end{center}
\end{figure}

A notional example of this refinement procedure is illustrated in Fig. \ref{fig:refinement}.
In this example, the set of candidate times used to compute $\boldsymbol\lambda_{est}$ is indicated by solid vertical lines.
It is evident that $g_{U(1,t)}(\mathbf\Gamma^T(t)\boldsymbol\lambda_{est}) \le 1$ at all of these times.
However, $g_{U(1,t)}(\mathbf\Gamma^T(t)\boldsymbol\lambda_{est}) \le 1-\epsilon_{remove}$ for two of these times (indicated by x).
These times are removed from $T^{est}$.
Next, the times of local maxima of $g_{U(1,t)}(\mathbf\Gamma^T(t)\boldsymbol\lambda_{est})$ that are greater than one (indicated by triangles) are added to $T^{est}$.
Because $\max_{t \in T}g_{U(1,t)}(\mathbf\Gamma^T(t)\boldsymbol\lambda_{est}) > 1+\epsilon_{cost}$, the solution is not within the specified tolerance, so further refinement is necessary.
Using the updated set of candidate times, $\boldsymbol\lambda_{est}$ is recomputed.
The evolution of $g_{U(1,t)}(\mathbf\Gamma^T(t)\boldsymbol\lambda_{est})$ for this updated estimate is shown as a dashed line.
It is evident that this new $\boldsymbol\lambda_{est}$ satisfies the convergence criteria, allowing the refinement to terminate.

A noteworthy advantage of this iteration procedure over other approaches such as Gilbert's algorithm \cite{Gilbert1971} is that each iteration constitutes a feasible solution with bounded sub-optimality.
Accordingly, this algorithm is well-suited for real-time applications with hard limits on computation time.

\subsection*{Extraction of Optimal Control Inputs}

Once a set of optimal control input times $T^{opt}$ and outward normal direction $\boldsymbol\lambda_{opt}$ are obtained, it is necessary to compute a set of optimal control inputs.
This can be accomplished by solving (\ref{eq:originalcontrolproblem}) considering control inputs only at times in $T_{opt}$.
However, this requires solving a convex optimization problem with complexity that depends on the selected cost function.

To minimize computation effort, the two-step approach described in the following simplifies the required optimization problem to a quadratic program.
First, an optimal control input direction is computed for each time in $T^{opt}$ for each control mode.
In the event that the optimal control input direction is not unique, the algorithm requires a set of points such that the convex hull includes all solutions of the support function (e.g. all columns of $\mathbf W$ that maximize the contact function for the examples in Table \ref{tab:examplecontactsupport}).
Second, the magnitudes of the control inputs are computed by solving a quadratic program that minimizes the error between $\boldsymbol w$ and the state reached by control inputs applied in the specified directions.
Provided that $\boldsymbol\lambda_{opt}$ is properly computed (i.e. the solver used in the iterative refinement algorithm converged), the residual error will be negligible for practical applications.
Additionally, the objective is formulated as the quadratic product of the error vector and a user-specified positive definite weight matrix $\mathbf Q$ to ensure well-behaved solutions.
This optimal control input extraction algorithm is described in the following pseudocode.\\

\noindent\textbf{Algorithm 3: Control Input Extraction}\\
\textbf{Inputs:} $T_{opt}$, $\boldsymbol\lambda_{opt}$, $\mathbf\Gamma(t)$, $\boldsymbol w$, $\mathbf Q$\\
\textbf{Outputs:} $\boldsymbol u_{opt}(t)$ \\
\textbf{loop} $t_j \in T^{opt}$\\
\indent $\hat{\boldsymbol u}_{opt, j}(t) \gets \delta(t-t_j)\boldsymbol s_{U(1,t_j)}(\mathbf\Gamma^T(t_j)\boldsymbol\lambda_{opt})$ \\
\indent $\boldsymbol y_j \gets \mathbf\Gamma(t_j)\boldsymbol s_{U(1,t_j)}(\mathbf\Gamma^T(t_j)\boldsymbol\lambda_{opt})$\\
$\boldsymbol\alpha \gets$ solution to optimization problem:\\
\indent minimize: $\boldsymbol w_{err} ^T\mathbf Q\boldsymbol w_{err} $\\
\indent subject to: $\boldsymbol w_{err} = \boldsymbol w - \sum\alpha_j\boldsymbol y_j, \: \alpha_j \ge 0, \: \sum \alpha_j \le \boldsymbol\lambda_{opt}^T\boldsymbol w$\\
$\boldsymbol{u}_{opt}(t) \gets \boldsymbol{0}$ \\
\textbf{loop} $t_j \in T^{opt}$\\
\indent $\boldsymbol u_{opt}(t) \gets \boldsymbol u_{opt}(t)+\alpha_j\hat{\boldsymbol u}_{opt, j}(t)$\\
\textbf{return} $\boldsymbol u_{opt}(t)$  \\

\begin{figure}[b]
\begin{center}
\includegraphics[trim={8cm 5cm 8cm 5cm},clip,width=0.8\columnwidth]{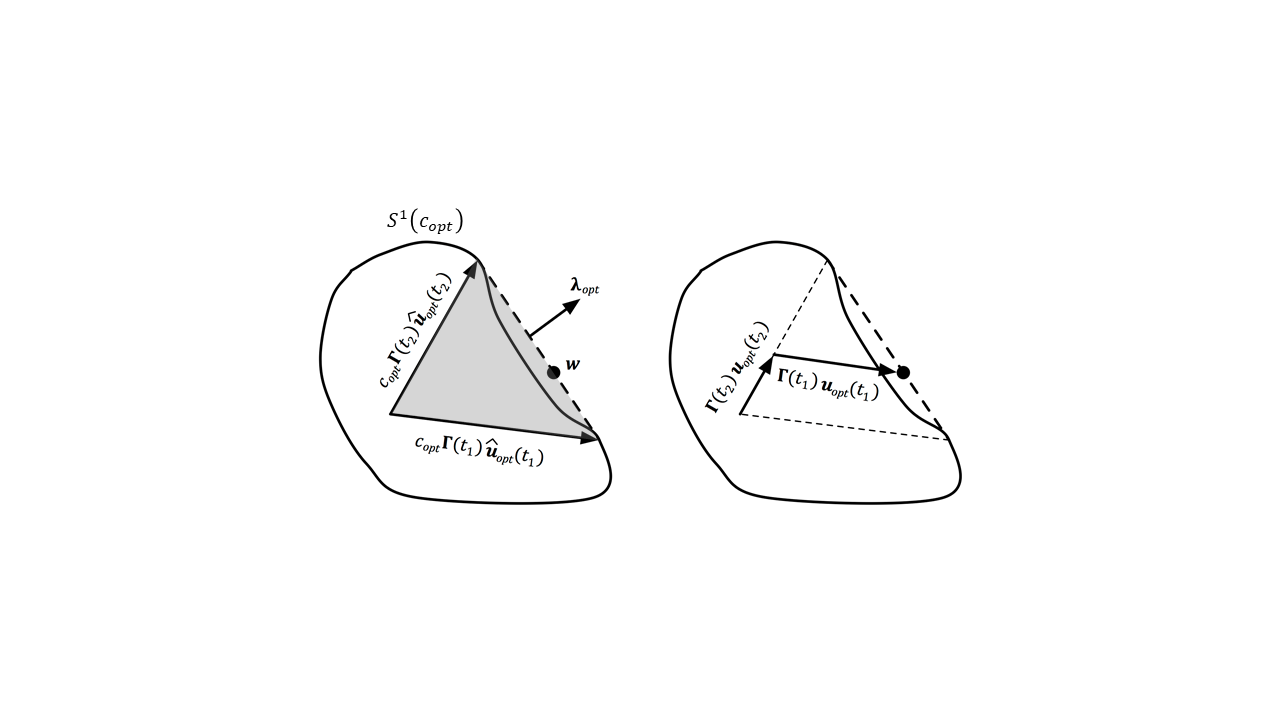}
\caption{Illustration of example optimal control input extraction for two-dimensional example including computation of optimal control input directions (left) and computation of control input magnitudes (right).}\label{fig:extraction}
\end{center}
\end{figure}

A notional example of the optimal control input extraction algorithm is shown in Fig. \ref{fig:extraction} for a two-dimensional system.
In this example, there are two candidate times for optimal control inputs.
For each of these times, the optimal control input direction $\hat{\boldsymbol u}_{opt}(t)$ is computed that reaches the supporting hyperplane at a cost of $c_{opt}$ as shown in the left plot.
The set of states that can be reached by a convex combination of these maneuvers is indicated by the shaded gray region.
Next, a convex combination of these control inputs is computed that reaches the specified $\boldsymbol w$ at cost $c_{opt}$.


\section{Validation and Performance Assessment}\label{sec:validation}

The proposed algorithm is validated through application to challenging spacecraft formation reconfiguration problems based on the Miniaturized Distributed Occulter/Telescope (mDOT) small satellite mission recently proposed by the authors \cite{Koenig2019, KoenigThesis}.
This mission uses a nanosatellite equipped with a telescope and a microsatellite equipped with a small occulter to obtain direct images of debris disks or large exoplanets from earth orbit. 
The technology demonstration variant of this mission requires autonomous formation reconfigurations that are challenging for three reasons: 1) the spacecraft have a large nominal separation of 500km established through a difference in the right ascension of the ascending node (RAAN), 2) the formation is deployed in a perturbed, eccentric orbit, and 3) the spacecraft are subject to time-varying attitude constraints to facilitate communication with ground stations.
The first two challenges are addressed by using the accurate linear dynamics model for $J_2$-perturbed relative motion provided in the appendix, which is based on mean absolute and relative orbital elements \cite{Koenig2017}.
Mean orbit elements are computed from osculating (i.e. instantaneous) orbit elements by applying a transformation that removes the effects of short-period perturbations.
This dynamics model is valid for arbitrarily large differences in right ascension of the ascending node (RAAN), enabling application on the mDOT mission.
The third challenge can be addressed by using the optimal impulsive control algorithm proposed in this paper.

\begin{figure}[tb]
\begin{center}
\includegraphics[trim={0cm 0cm 0cm 0cm},clip,width=0.9\columnwidth]{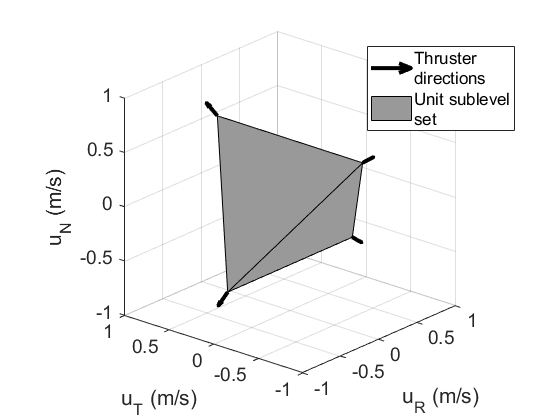}
\caption{Illustration of thruster orientations in RTN frame and set of admissible maneuvers with unit cost in fixed-attitude mode.}\label{fig:tetra}
\end{center}
\end{figure}

To validate the proposed algorithm, it is first necessary to define the reconfiguration problem.
It is assumed that the reconfigurations start at the apogee of an orbit with a 25000km semimajor axis, resulting in an orbit period of 10.92 hours.
Formation reconfigurations are allowed to take three orbits.
Thus, $t_i$ is selected as 0 and $t_f$ is selected as 117990 seconds.
The control domain $T$ is selected as a uniform discretization of the interval $[t_i, t_f]$ with thirty second intervals for a total of 3934 candidate control input times.
The occulter spacecraft (which performs all maneuvers) will need to communicate with a ground station every orbit to downlink data from each observation.
To accommodate this constraint, it is assumed that the occulter must maintain a fixed attitude in the radial/tangential/normal (RTN) frame for a two-hour window surrounding the perigee of each orbit to facilitate communications with ground stations.
While this interval is significantly longer than typical ground contacts, this choice helps to illustrate the different behavior of $g_{U(1,t)}(\mathbf\Gamma^T(t)\boldsymbol\lambda)$ in each mode (with and without attitude constraints).
Additionally, it is assumed that the occulter has four thrusters arranged in an equilateral tetrahedral configuration.
The alignment of each of these thrusters in the RTN frame in the fixed-attitude mode are the columns of the matrix $\mathbf{V}^{vertex}$ defined as
\begin{equation}
\label{eq:thrusterlocations}
\mathbf{V}^{vertex} = 
\begin{bmatrix} 
\sqrt{2/3} & -\sqrt{2/3} & 0 & 0 \\
0 & 0 & \sqrt{2/3}  & -\sqrt{2/3} \\
-\sqrt{1/3} & -\sqrt{1/3} & \sqrt{1/3}  & \sqrt{1/3}  
\end{bmatrix}
\end{equation}
The corresponding $\boldsymbol V^{face}$ used to define the cost function is given by
\begin{equation}
\label{eq:mdotsublevel}
\mathbf{V}^{face} = 
\frac{1}{3}\begin{bmatrix} 
-\sqrt{2/3} & 0 & \sqrt{1/3} \\
\sqrt{2/3} & 0 & \sqrt{1/3} \\
0 & -\sqrt{2/3} & -\sqrt{1/3} \\
0 & \sqrt{2/3} & -\sqrt{1/3} \\
\end{bmatrix}
\end{equation}
The thruster directions and the set of admissible maneuvers of unit cost for the fixed attitude mode are illustrated in Fig. \ref{fig:tetra}.
The row vectors of $\mathbf V^{face}$ are parallel to the outward normal directions to the faces of the tetrahedron in this plot.
It is also assumed that no attitude constraints are enforced outside of the two-hour window surrounding the orbit perigee.
During these times, it assumed that the occulter spacecraft can freely rotate so that any maneuver can be executed by firing a single thruster.

Under these assumptions, the time domains for the fixed-attitude mode (denoted $T_1$) and unconstrained mode (denoted $T_2$) are given by
\begin{align*}
\begin{split}
    T_{1} = \{ t : t \in T, \; |t-(N+0.5)T_{orbit}| < 1\textrm{hr}, \; N \in \mathbb{Z} \} \\
    T_{2} = \{ t : t \in T, \; |t-(N+0.5)T_{orbit}| \ge 1\textrm{hr}, \; \forall N \in \mathbb{Z} \}
\end{split}
\end{align*}
and $T_{orbit}$ is the orbit period.
The cost function for this reconfiguration problem is given by
\begin{equation}
\label{eq:mdotgdef}
f(\boldsymbol u, t) = 
\begin{Bmatrix}
\max(\mathbf V^{face}\boldsymbol u) & t \in T_1 \\
||\boldsymbol u||_2 & t \in T_2 
\end{Bmatrix}
\end{equation}
The contact and support functions for $U(1,t)$ take the forms provided in Table \ref{tab:examplecontactsupport} for the corresponding form of the cost function at time $t$.

Key algorithm parameters are described in the following.
For the initialization algorithm, the provided $T^d$ includes 20 times evenly distributed between $t_i$ and $t_f$.
The provided $\boldsymbol\lambda_{est}$ is a unit vector parallel to the target pseudostate.
The initial set of candidate times $T^{est}$ includes the six times in $T_d$ at which $g_{ U(1,t)}(\mathbf\Gamma^T(t)\boldsymbol\lambda_{est})$ is largest.
The tolerances $\epsilon_{cost}$ and $\epsilon_{remove}$ in the refinement algorithm were selected as 0.01.
Finally, the error weight matrix $\mathbf Q$ is the identity matrix.
The algorithms defined in the previous section were implemented in MATLAB and CVX was used to solve the required convex optimization problems in the iterative refinement and optimal control input extraction algorithms \cite{cvx,gb08}.

%
%
%
%
%
%
In all test cases presented in the following, the normalized residual error ($||\boldsymbol w_{err}||_2/||\boldsymbol w||_2$) was less than 0.01\%, indicating that the solver reliably converged for both the iterative refinement algorithm and the maneuver extraction algorithm.


\subsection*{Example Formation Reconfiguration Problem}

The proposed algorithm is first used to compute an optimal maneuver sequence for an example formation reconfiguration problem over three orbits.
The initial absolute orbit (used to evaluate $\boldsymbol\Gamma(t)$ using the dynamics model in the appendix) and target pseudostate are provided in Table \ref{tab:exampleproblem}.
The target pseudostate is scaled by the orbit semimajor axis in the table so that it is numerically comparable to the change in the relative position.
\begin{table}[htb]
\centering
\caption{Initial mean orbit and target pseudostate.}
\begin{tabular}{c c c c c c}
\hline
\multicolumn{6}{c}{Initial mean absolute orbit \textbf{\oe}$_c(t_i)$} \\
\multicolumn{1}{c}{$a(km)$} & \multicolumn{1}{c}{$e(-)$} & \multicolumn{1}{c}{$i(^o)$} & \multicolumn{1}{c}{$\Omega(^o)$} & \multicolumn{1}{c}{$\omega(^o)$} & \multicolumn{1}{c}{$M(^o)$} \\ \hline
25000   & 0.7 & 40 & 358  &  0  & 180    \\  \hline
\multicolumn{6}{c}{Target pseudostate $\boldsymbol w$ (m)} \\
\multicolumn{1}{c}{$a\delta a$} & \multicolumn{1}{c}{$a\delta\lambda$} & \multicolumn{1}{c}{$a\delta e_x$} & \multicolumn{1}{c}{$a\delta e_y$} & \multicolumn{1}{c}{$a\delta i_x$} & \multicolumn{1}{c}{$a\delta i_y$} \\ \hline
50 & 5000 & 100 & 100 & 0 & 400   \\ \hline
\end{tabular}\label{tab:exampleproblem}
\end{table}

A solution that reaches the target pseudostate and satisfies the optimality criteria to within a tolerance of $\epsilon_{cost}$ was found using only three iterations of Algorithm 2.
The optimal dual variable is given by  $\boldsymbol\lambda_{opt} = 10^{-6}\times[34.97, 3.42, 30.68, 17.84, -9.34, 146.79]^T$.
%
%
The delta-v cost of the computed maneuver sequence is 82.4mm/s, which is within the specified 1\% tolerance of the lower bound of 82.0mm/s computed by evaluating $\boldsymbol\lambda_{opt}^T\boldsymbol w/g_{S(1)}(\boldsymbol\lambda_{opt})$.
The optimal maneuver sequence consists of the three maneuvers in the RTN frame provided in Table \ref{tab:optexamplemaneuvers}.
It is noteworthy that these maneuvers include significant radial components, which contradicts the expected behavior from the closed-form solutions developed by Chernick \cite{Chernick2018}.
This behavior arises from the fact that Chernick's solutions assume that in-plane ($\delta a$, $\delta\lambda$, $\delta e_x$, and $\delta e_y$) and out-of-plane ($\delta i_x$ and $\delta i_y$) control are decoupled, while this algorithm optimally couples in-plane and out-of-plane control.
\begin{table}[htb]
\centering
\caption{Optimal maneuvers for example scenario.}
\begin{tabular}{ c  r  r  r }
\hline
$t_j$ (sec)                        &    16050   &    23280   &   107100\\ \hline
$u_R(t_j)$ (mm/s)                      & 9.68 &   0.00 &  16.51 \\
$u_T(t_j)$ (mm/s)                      & -23.02 & -0.40  & 15.68  \\
$u_N(t_j)$ (mm/s)                      & -25.56 &  -0.04 &  40.26  \\ \hline
\end{tabular}\label{tab:optexamplemaneuvers}
\end{table}

\begin{figure}[htb]
\begin{center}
\includegraphics[trim={6.5cm 4cm 8cm 4cm},clip,width=0.9\columnwidth]{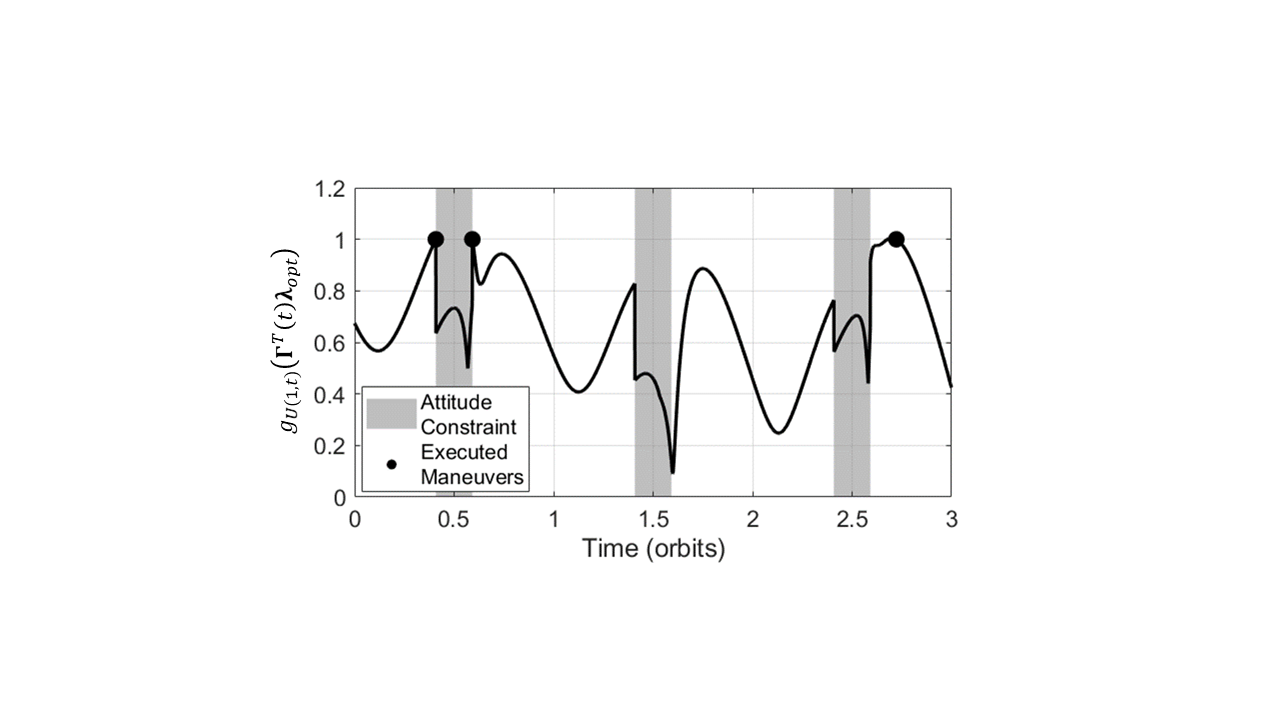}
\caption{Evolution of $g_{U(1,t)}(\mathbf\Gamma(t)\boldsymbol\lambda_{opt})$ for example problem including optimal maneuver times (black circles) and attitude constraints (gray).}\label{fig:primerexample}
\end{center}
\end{figure}

The evolution of $g_{U(1,t)}(\mathbf\Gamma(t)\boldsymbol\lambda_{opt})$ for this solution is illustrated in Fig. \ref{fig:primerexample}.
The optimal maneuver times are indicated by black circles and the time intervals in which the fixed attitude constraint is enforced are indicated by gray shading.
It is evident from this plot that the optimality criteria are satisfied to within the specified tolerance of 0.01.


\subsection*{Comparison to State-of-the-Art Approaches}

To demonstrate the improved computational efficiency of the proposed algorithm, Monte Carlo simulations were conducted using the proposed algorithm and two other state-of-the art algorithms.
The first reference algorithm is direct optimization (i.e. solving the problem in (\ref{eq:primalprob})) where control inputs are allowed at all 3934 times in the discretization of $T$.
This requires solving a single convex optimization problem with 11802 variables and 6 equality constraints.
The second reference algorithm is an implementation Gilbert's algorithm \cite{Gilbert1971} that has been modified to use the support and contact functions derived in Section \ref{sec:contactsupport}.
This algorithm is selected because it is the fastest algorithm in literature that provides guaranteed convergence and is applicable to this class of problem.
Gilbert's algorithm is applied with an error tolerance of $0.1$m in the reached pseudostate and an overexpansion factor of 0.01 to reduce the number of required iterations.
To provide a fair comparison between the run times, CVX \cite{cvx,gb08} is used to solve the quadratic programs used in Gilbert's refinement procedure.

Each algorithm was used to compute optimal impulsive control inputs for 200 pseudostates randomly selected from a zero-mean Gaussian distribution with standard deviations of 1km for each ROE.
The minimum, mean, and maximum computation time and number of required iterations for each of these algorithms run on a desktop computer with a 3.5 GHz processor are provided in Table \ref{tab:runtime}.
\begin{table}[b]
\centering
\caption{Minimum, mean, and maximum values for computation time and required iterations in Monte Carlo simulations.}
\begin{tabular}{ l c c c c c c }
\hline
& \multicolumn{3}{c}{Computation time (s)} & \multicolumn{3}{c}{\# Iterations} \\
Algorithm & Min & Mean & Max  & Min & Mean & Max \\ \hline
Proposed  & 0.76 & 2.43 & 4.32 & 1 & 3.99 & 8 \\ 
Gilbert   & 3.72 & 7.36 & 14.7 & 9 & 19.1 & 40\\
Direct    & 85.5 & 90.9 & 98.5 & - & - & -\\ \hline
\end{tabular}\label{tab:runtime}
\end{table}
It is clear from these tests that both the algorithm proposed in this paper and Gilbert's algorithm are more than an order of magnitude faster than direct optimization.
While the proposed algorithm required approximately one third of the computation time of Gilbert's algorithm, it should be noted that this behavior may change if different solvers are chosen.
However, it is clear that the algorithm proposed in this paper requires on average five times fewer iterations than Gilbert's algorithm, and correspondingly fewer evaluations of the contact function.
It is event from this comparison that the proposed algorithm is well-suited to problems where evaluation of the contact and support functions is computationally expensive.


\subsection*{Sensitivity to Poor Initialization}

\begin{figure}[htb]
\begin{center}
\includegraphics[trim={1cm 0cm 1.5cm 0cm},clip,width=1.0\columnwidth]{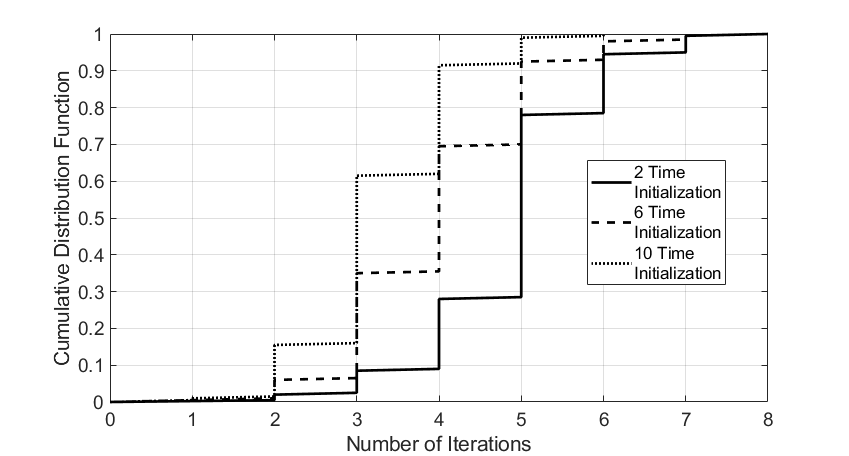}
\caption{Cumulative distribution functions for the number of required refinement iterations in Monte Carlo simulations for three initialization schemes.}\label{fig:numiters}
\end{center}
\end{figure}

To characterize the sensitivity of the computation cost to the initial set of candidate times, the Monte Carlo simulations using the proposed algorithm were repeated with two different sets of initial candidate control input times.
The first initial set of times includes only $t_i$ and $t_f$.
This initialization is intended to capture the worst-case computation cost because it is unlikely that the optimal cost can be reached with only two maneuvers.
The second initialization includes ten candidate times evenly spaced in the interval $[t_i, t_f]$.
This initialization ensures that the initial candidate times are reasonably close to optimal times, but requires the algorithm to check a larger number of constraints in the iterations.
The initializations with two, six, and ten candidate times required averages of 4.90, 3.99, and 3.31 iterations of Algorithm 2, respectively.
Fig. \ref{fig:numiters} shows the distribution of the number of iterations required to solve these reconfiguration problems for all three initialization schemes.
It is evident that increasing the number of times in the initial set of candidate control input times provides a modest decrease in the number of required iterations.
However, the algorithm reliably converges from a worst-case initialization in six iterations or less for most problems.
Also, including more candidate times increases the complexity of the optimization problems that must be solved in each iteration.
Thus, the ideal number of candidate times for initialization will depend on the limitations of available solvers for a specified application.
Overall, these results show that the algorithm is robust to poor initialization data with a nominal increase of approximately 30\% in the number of required iterations.


\subsection*{Sensitivity to Time Discretization}

\begin{figure}[htb]
\begin{center}
\includegraphics[trim={0cm 0cm 0cm 0},clip,width=\columnwidth]{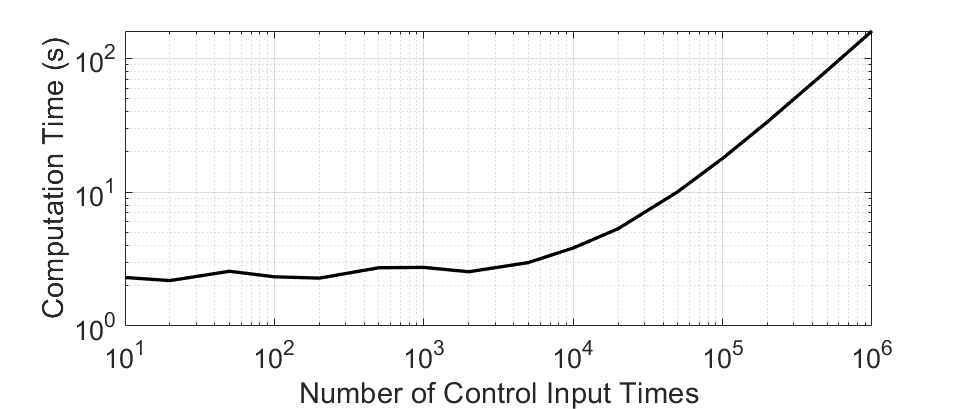}
\caption{Computation time vs. number of candidate control input times.}\label{fig:timevstimes}
\end{center}
\end{figure}

%
To characterize the relationship between the computation time and the density of the time discretization, a selection of ten problems from the Monte Carlo experiment were solved with the the number of times in $T$ varied between ten and one million.
The maximum computation time of the ten problems is plotted against the number of times in $T$ in Fig. \ref{fig:timevstimes}.
It is evident from this plot that the required computation time is nearly constant if the number of times is $10^4$ or less.
Beyond this level, the cost increases linearly with the number of samples.
This is because the majority of the computation effort is spent updating the set of candidate optimal control input times, which has equivalent computation cost to one evaluation of the contact function for the reachable set.
Since the computation time for direct optimization algorithms varies exponentially with the number of variables, it is evident that the algorithm proposed in this paper is well suited to problems that require fine discretizations of the time domain.

These simulations demonstrate that the proposed algorithm exhibits robust convergence to globally optimal control inputs several times faster than comparable algorithms in literature.
The authors have also found in preliminary tests that the proposed algorithm runs on a 800MHz CubeSat-compatible microprocessor with run times of under ten seconds, making it suitable for deployment in real-time applications such as the guidance, navigation, and control system for the VIrtual Super-resolution Optics with Reconfigurable Swarms (VISORS) mission \cite{Koenig2021}.


\section{Conclusions}\label{sec:conclusions}

This paper proposes a new simple and robust solution methodology for a class of fuel-optimal impulsive control problems for linear systems with time-varying cost.  
First, the properties of the reachable set are derived.
Next, contact and support functions for the reachable set are derived, enabling use of existing algorithms for optimization over parameterized sets.
After reformulating the problem as a semi-infinite convex program, it is shown that this class of problem is identical to the form studied in primer vector theory under the assumption that the cost is a constant $p$-norm.
It follows that the contact function is a generalization of the norm of the primer vector and provides a geometric interpretation of the time-invariant portion of the primer vector.
Finally, a three-step algorithm is proposed that provides efficient and robust computation of globally fuel-optimal impulsive control input sequences.
The geometry of the problem is leveraged in every step to reduce computational cost and ensure robustness.

The algorithm is validated through implementation in challenging spacecraft formation reconfiguration problems based on the proposed Miniaturized Distributed Occulter/Telescope small satellite mission.
It is found that the algorithm is able to compute a maneuver sequence with a total cost within 1\% of the global optimum within eight iterations in all test cases, including those with worst-case initializations.
Also, the normalized residual error of all computed solutions was no larger than 0.01\%, indicating reliable convergence.
It was found that the algorithm is more than an order of magnitude faster than direct optimization and three times faster than the best indirect optimization algorithm available in literature.
Indeed, the required computation time is nearly constant unless the discretization of the time domain includes more than ten thousand samples.

Overall, the proposed solution methodology provides a real-time-capable means of computing globally optimal impulsive control input sequences for a wide range of linear time-variant dynamical systems with complex cost behaviors.
%
%

\section*{Acknowledgment}

This work was supported by the Air Force Research Laboratory’s Control, Navigation, and Guidance
for Autonomous Systems (CoNGAS) project under contract FA9453-16-C-0029 and by a NASA Office of the Chief Technologist Space Technology Research Fellowship, NASA Grant No. NNX15AP70H.

\appendix

\section{Dynamics Model}\label{sec:dynamicsmodel}

This appendix describes the dynamics model for control of spacecraft formations in $J_2$-perturbed orbits that is used in Section \ref{sec:validation}.
This model is valid for mean orbit elements, which are computed from the osculating (or instantaneous) orbit elements by applying a transformation that removes short-period oscillations.
Similar models have been used for spacecraft formation control on multiple missions such as TanDEM-X and PRISMA to avoid excess fuel expenditure due to short-period perturbations \cite{DAmicoThesis, Gaias2015diffdrag1}.

Let $\mu$ denote earth's gravitational parameter, $R_E$ denote earth's mean radius, and $J_2$ denote earth's second degree zonal geopotential coefficient.
These constants are given by
\begin{align*}
\begin{split}
\mu = 3.986 \times 10^{14}\; \textrm{m}^3/\textrm{s}^2, \qquad R_E = 6.378 \times 10^{6} \; \textrm{m}, \\
J_2 = 1.082 \times 10^{-3}. \qquad \qquad \quad
\end{split}
\end{align*}

The mean absolute orbits for each spacecraft are described by Keplerian orbit elements, which include the semimajor axis $a$, eccentricity $e$, inclination $i$, right ascension of the ascending node $\Omega$, argument of perigee $\omega$, and mean anomaly $M$.
Using these elements, the orbit element vector \textbf{\oe} is defined as
\begin{align*}
\textbf{\oe} = \begin{bmatrix} a & e & i & \Omega & \omega & M \end{bmatrix}^T.
\end{align*}
The unforced dynamics for the mean orbit elements of a spacecraft in a $J_2$-perturbed orbit are given in \cite{Battin1987} as
\begin{equation}
\label{eq:j2dynamics}
\begin{split}
\dot{\textbf{\oe}} = 
\begin{pmatrix} 
\dot{a} \\ \dot{e} \\ \dot{i} \\ \dot{\Omega} \\ \dot{\omega} \\ \dot{M}
\end{pmatrix}=
\begin{pmatrix}
0 \\ 0 \\ 0 \\
-\frac{3 J_2 R_E^2\sqrt{\mu}}{2 a^{7/2}\eta^4} \cos(i)  \\
\frac{3 J_2 R_E^2\sqrt{\mu}}{4 a^{7/2}\eta^4} (5\cos^2(i) -1) \\
\sqrt{\frac{\mu}{a^3}}+\frac{3 J_2 R_E^2\sqrt{\mu}}{4 a^{7/2}\eta^3} (3\cos^2(i) -1)
\end{pmatrix}.
\end{split}
\end{equation}
It is evident from this equation that $a$, $e$, and $i$ are constant while $\Omega$, $\omega$, and $M$ vary linearly with time.
Accordingly, the absolute orbit at time $t$ is related to the orbit at time $t_i$ by
\begin{align*}
    \textbf\oe(t) = \textbf\oe(t) + (t-t_i)\dot{\textbf{\oe}}
\end{align*}

The mean relative orbital elements (ROE) state used in this paper is defined with respect to the mean orbits of the chief, denoted by subscript $c$, and the deputy, denoted by subscript $d$, by
\begin{equation}
\label{eq:roedef}
\begin{split}
\boldsymbol x =
\begin{pmatrix} \delta a \\ \delta\lambda \\ \delta e_x \\ \delta e_y \\ \delta i_x \\ \delta i_y \end{pmatrix} = 
\begin{pmatrix} \Delta a/a_c \\ \Delta M + \eta_c(\Delta\omega + \Delta\Omega\cos (i_c)) \\ e_d\cos(\omega_d) - e_c\cos(\omega_c) \\ e_d\sin(\omega_d) - e_c\sin(\omega_d) \\ \Delta i \\ \Delta\Omega\sin (i_c) \end{pmatrix}
\end{split}
\end{equation}
where $\eta = \sqrt{1-e^2}$ and the operator $\Delta$ denotes the difference between the orbit elements of the deputy and chief (e.g. $\Delta a = a_d-a_c$).
Without loss of generality, it is assumed that the telescope spacecraft is the chief and the occulter spacecraft is the deputy and all maneuvers are executed by the occulter.

The dynamics model includes a control input matrix $\mathbf B(t)$ and a state transition matrix (STM) $\mathbf\Phi(t, t_f)$.
The control input matrix used in this paper relates the effect of an applied impulse to its effect on the ROE, which is given in \cite{Chernick2018} as
\begin{align*}
\mathbf B(t) = \sqrt{\frac{a_c}{\mu}}\begin{bmatrix} 
B_{11} & B_{12} & 0 \\
B_{21} & 0 & 0 \\
B_{31} & B_{32} & B_{33} \\
B_{41} & B_{42} & B_{43} \\
0 & 0 & B_{53} \\
0 & 0 & B_{63}
\end{bmatrix}.
\end{align*}
The nonzero terms of this matrix are given by
\begin{align*}
\begin{split}
B_{11} = \frac{2}{\eta}e\sin(\nu), \quad
B_{12} = \frac{2}{\eta}(1+e\cos(\nu)), \\
B_{21} = -\frac{2\eta^2}{1+e\cos(\nu)}, \quad
B_{31} = \eta\sin(\theta), \\
B_{32} = \eta\frac{(2+e\cos(\nu))\cos(\theta)+e\cos(\omega)}{1+e\cos(\nu)}, \\
B_{33} = \frac{\eta e\sin(\omega)\sin(\theta)}{\tan (i) (1+e\cos(\nu))}, \quad B_{41} = -\eta\cos(\theta), \\
B_{42} = \eta\frac{(2+e\cos(\nu))\sin(\theta)+e\sin(\omega)}{1+e\cos(\nu)}, \\
B_{43} = -\frac{\eta e\cos(\omega)\sin(\theta)}{\tan (i) (1+e\cos(\nu))}, 
\end{split}
\end{align*}
\begin{align*} 
\begin{split}
B_{53} = \frac{\eta\cos(\theta)}{1+e\cos(\nu)}, \quad
B_{63} = \frac{\eta\sin(\theta)}{1+e\cos(\nu)}
\end{split}
\end{align*}
where $\theta = \omega+\nu$ and $\nu$ is the true anomaly, which is related to the mean anomaly by Kepler's equation.
The columns of the control matrix correspond to thrusts applied to the deputy spacecraft in the radial (R), along-track (T), and cross-track (N) directions, respectively.
The R direction is aligned with the position vector of the spacecraft, the N direction is aligned with the angular momentum vector of the orbit, and the T direction completes the right-handed triad.
All orbit elements are evaluated at the time at which the impulse is applied.

As demonstrated in \cite{Koenig2017}, the STM for the state defined in (\ref{eq:roedef}) can be computed by performing a first order Taylor expansion on the equations of relative motion and solving the linearized system of equations in closed-form.
For the ROE definition in (\ref{eq:roedef}) and dynamics model in (\ref{eq:j2dynamics}), the resulting STM is given by
\begin{align*}
\begin{split}
\mathbf\Phi(t, t_f) = 
\begin{bmatrix} 
\Phi_{11} & 0 & 0 & 0 & 0 & 0 \\
\Phi_{21} & \Phi_{22} & \Phi_{23} & \Phi_{24} & \Phi_{25} & 0 \\
\Phi_{31} & 0 & \Phi_{33} & \Phi_{34} & \Phi_{35} & 0 \\
\Phi_{41} & 0 & \Phi_{43} & \Phi_{44} & \Phi_{45} & 0 \\
0 & 0 & 0 & 0 & \Phi_{55} & 0 \\
\Phi_{61} & 0 & \Phi_{63} & \Phi_{64} & \Phi_{65} & \Phi_{66}
\end{bmatrix}.
\end{split}
\end{align*}
The nonzero terms of this STM are given by
\begin{align*}
\begin{split}
\Phi_{11} = 1, \quad 
\Phi_{21} =(-1.5\sqrt{\mu/a^3}\Delta t-7\kappa\eta P) \Delta t, \\
\Phi_{22} = 1, \quad
\Phi_{23} = 7\kappa e_{x1} P\Delta t/\eta, \quad
\Phi_{24} =  7\kappa e_{y1} P\Delta t /\eta, \\
\Phi_{25} = -7\kappa\eta S \Delta t, \quad
\Phi_{31} = 3.5\kappa e_{y2} Q\Delta t, \\
\Phi_{33} = \cos(\dot{\omega}\Delta t)-4\kappa e_{x1}e_{y2}GQ\Delta t, \\
\Phi_{34} = -\sin(\dot{\omega}\Delta t)-4\kappa e_{y1}e_{y2}GQ\Delta t, \\
\Phi_{35} = 5\kappa e_{y2} S \Delta t, \quad
\Phi_{41} = -3.5\kappa e_{x2} Q\Delta t, \\
\Phi_{43} = \sin(\dot{\omega}\Delta t)+4\kappa e_{x1}e_{x2}GQ\Delta t, \\
\Phi_{44} =  \cos(\dot{\omega}\Delta t)+4\kappa e_{y1}e_{x2}GQ\Delta t, \\
\Phi_{45} = -5\kappa e_{x2}S\Delta t, \quad
\Phi_{55} = 1, \quad
\Phi_{61} =  3.5\kappa S \Delta t, \\
\Phi_{63} =  -4\kappa e_{x1} GS\Delta t, \quad
\Phi_{64} =  -4\kappa e_{y1} GS\Delta t, \\
\Phi_{65} =   2\kappa T\Delta t, \quad
\Phi_{66} =  1
\end{split}
\end{align*}
where the substitutions given by
\begin{align*}
\begin{split}
\Delta t = t_f-t, \quad
\kappa = \frac{3J_2 R_E^2 \sqrt{\mu}}{4a^{7/2}\eta^4},  \quad 
G = \eta^{-2}, \\
P = 3\cos^2(i)-1, \quad
Q = 5\cos^2(i)-1, \quad
S = \sin(2i), \\ 
T = \sin^2(i), \quad
e_{x1} = e\cos(\omega(t)), \quad
e_{y1} = e\sin(\omega(t)), \\
e_{x2} = e\cos(\omega(t_f)), \quad
e_{y2} = e\sin(\omega(t_f))
\end{split}
\end{align*}
are used to simplify notation.
This STM is the same as the $J_2$-perturbed STM developed in \cite{Koenig2017} for quasi-nonsingular ROE except that the second row is modified to accommodate the changed definition of $\delta\lambda$.

\bibliographystyle{ieeetr}
\bibliography{root}

\begin{IEEEbiography}[{\includegraphics[width=1in,height=1.25in,clip,keepaspectratio]{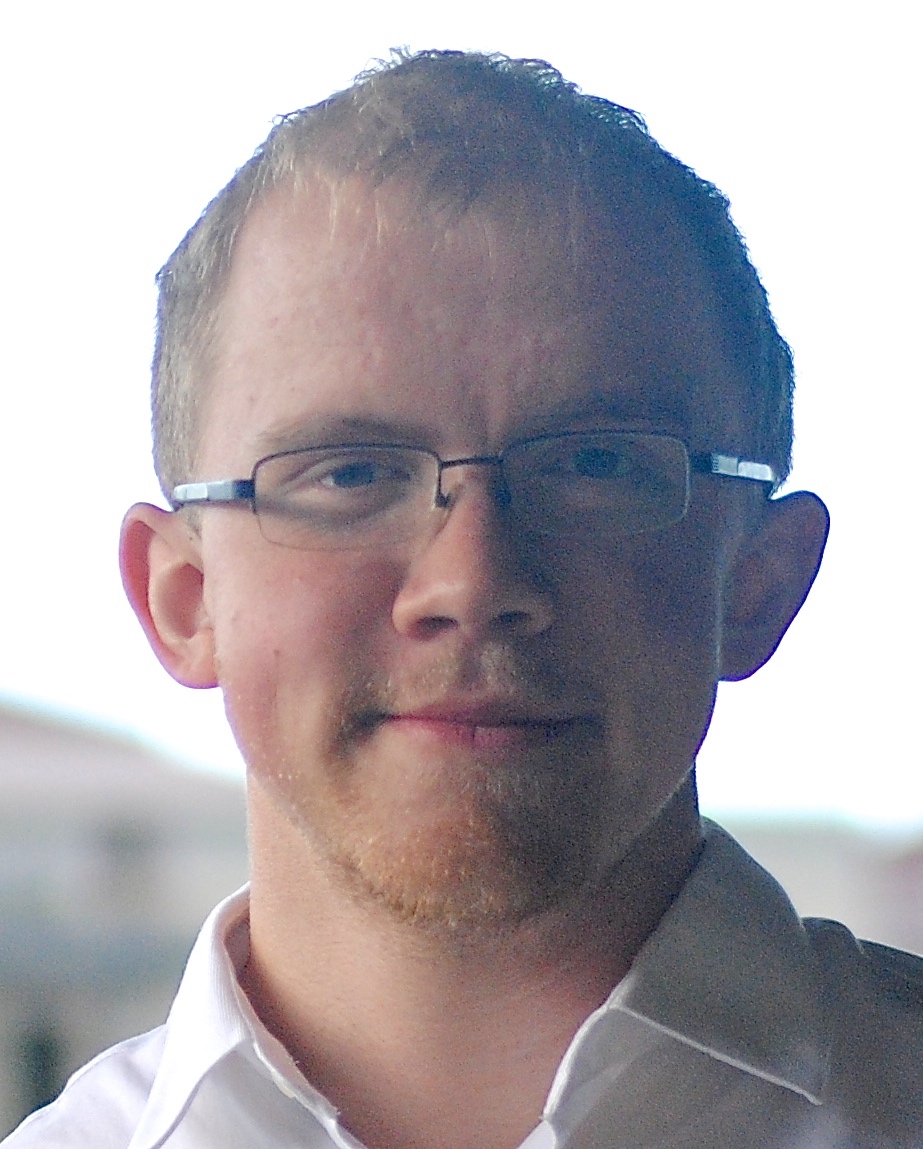}}]
{Adam W. Koenig} received his B.S. degree in aerospace engineering from Wichita State University in 2012 and his M.S. (2015) and Ph.D. (2019) degrees in aeronautics and astronautics from Stanford University.

Dr. Koenig is currently a postdoctoral scholar in the Space Rendezvous Laboratory at Stanford University.
His research interests include astrodynamics and advanced guidance, navigation, and control algorithms for distributed space systems.
His current activities include development of an angles-only navigation software payload for the StarFOX experiment on the Starling1 mission in development at NASA Ames Research Center, development of the guidance, navigation, and control algorithms for the VIrtual Super-resolution Optics with Reconfigurable Swarms (VISORS) mission funded by NSF, and design and analysis of the Miniaturized Distributed Occulter/Telescope (mDOT) mission concept for direct imaging of debris disks and large planets orbiting nearby stars.
\end{IEEEbiography}

\begin{IEEEbiography}[{\includegraphics[width=1in,height=1.25in,clip,keepaspectratio]{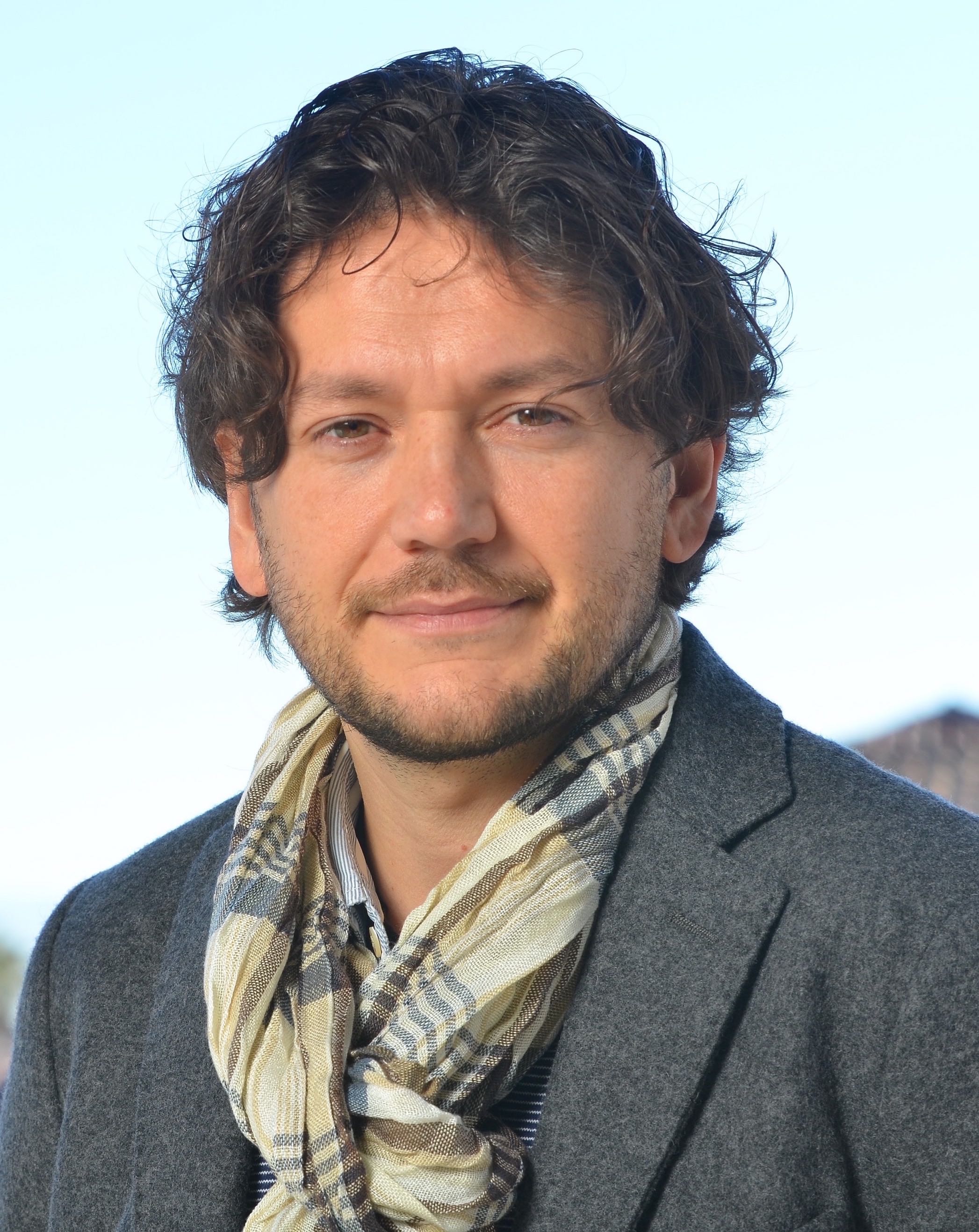}}]
{Simone D'Amico} received the B.S. and M.S. degrees from Politecnico di Milano (2003) and the Ph.D. degree from Delft University of Technology (2010). 
From 2003 to 2014, he was research scientist and team leader at the German Aerospace Center (DLR). 
There, he gave key contributions to the design, development, and operations of spacecraft formation-flying and rendezvous missions such as GRACE (United States/Germany), TanDEM-X (Germany), PRISMA (Sweden/Germany/France), and PROBA-3 (ESA). 

Since 2014, he has been Assistant Professor of Aeronautics and Astronautics at Stanford University, founding director of the Space Rendezvous Laboratory (SLAB), and Satellite Advisor of the Student Space Initiative (SSSI), Stanford’s largest undergraduate organization. 
He has over 150 scientific publications and 2000 Google Scholar citations, including conference proceedings, peer-reviewed journal articles, and book chapters. 
D'Amico's research aims at enabling future miniature distributed space systems for unprecedented science and exploration. 
His efforts lie at the intersection of advanced astrodynamics, GN\&C, and space system engineering to meet the tight requirements posed by these novel space architectures. 
The most recent mission concepts developed by Dr. D'Amico are a miniaturized distributed occulter/telescope (mDOT) system for direct imaging of exozodiacal dust and exoplanets and the Autonomous Nanosatellite Swarming (ANS) mission for characterization of small celestial bodies.
D’Amico’s research is supported by NASA, AFRL, AFOSR, KACST, and Industry. 
%
%
He is member of the Space-Flight Mechanics Technical Committee of the AAS, Associate Fellow of AIAA, and Associate Editor of the Journal of Guidance, Control, and Dynamics. 
%
\end{IEEEbiography}

\end{document}